%% file: RBIN.tex
\documentclass[reqno,twoside]{birkmult}

\usepackage{general,resistance,ifthen}
\usepackage{cancel,xy}

\usepackage[bookmarks, colorlinks=true, pdfstartview=FitV, linkcolor=blue, citecolor=blue, urlcolor=blue]{hyperref}
\usepackage{cite}  

\newcommand{\clean}{true}  
\newcommand{\version}[2]{\ifthenelse{\equal{\clean}{true}}{#1}{{\footnotesize #2}}}
\newcommand{\linenopax}{}

\numberwithin{equation}{section} \numberwithin{theorem}{section}

\begin{document}

\title{Resistance boundaries of infinite networks}

\author[P. E. T. Jorgensen]{Palle E. T. Jorgensen}
\address{University of Iowa, Iowa City, IA 52246-1419 USA}
\email{jorgen@math.uiowa.edu}

\author[E. P. J. Pearse]{Erin P. J. Pearse}
\address{University of Iowa, Iowa City, IA 52246-1419 USA}
\email{epearse@math.uiowa.edu}

\thanks{The work of PETJ was partially supported by NSF grant DMS-0457581. The work of EPJP was partially supported by the University of Iowa Department of Mathematics NSF VIGRE grant DMS-0602242.}

\begin{abstract}
A resistance network is a connected graph $(G,c)$. The conductance function $c_{xy}$ weights the edges, which are then interpreted as conductors of possibly varying strengths. The Dirichlet energy form $\mathcal E$ produces a Hilbert space structure ${\mathcal H}_{\mathcal E}$ on the space of functions of finite energy.  

The relationship between the natural Dirichlet form $\mathcal E$ and the discrete Laplace operator $\Delta$ on a finite network is given by $\mathcal E(u,v) = \la u, \Lap v\ra_2$, where the latter is the usual $\ell^2$ inner product. 
We describe a reproducing kernel $\{v_x\}$ for $\mathcal E$ and use it to extend the discrete Gauss-Green identity to infinite networks:
\[{\mathcal E}(u,v) = \sum_{G} u \Delta v + \sum_{\operatorname{bd}G} u \tfrac{\partial v}{\partial \mathbf{n}} ,\]
where the latter sum is understood in a limiting sense, analogous to a Riemann sum. This formula yields a boundary sum representation for the harmonic functions of finite energy. 

Techniques from stochastic integration allow one to make the boundary $\operatorname{bd}G$ precise as a measure space, and give a boundary integral representation (in a sense analogous to that of Poisson or Martin boundary theory). This is done in terms of a Gel'fand triple $S \ci {\mathcal H}_{\mathcal E} \ci S'$ and gives a probability measure $\mathbb{P}$ and an isometric embedding of ${\mathcal H}_{\mathcal E}$ into $L^2(S',\mathbb{P})$, and yields a concrete representation of the boundary as a set of linear functionals on $S$.

\end{abstract}

  \keywords{Dirichlet form, graph energy, discrete potential theory, graph Laplacian, weighted graph, tree, electrical resistance network, effective resistance, resistance form, Markov process, random walk, transience, Martin boundary, boundary theory, boundary representation, harmonic analysis, Hilbert space, orthogonality, unbounded linear operators, reproducing kernel. \\ {\footnotesize Date: \bf\today}\vstr}

  \subjclass{
    Primary:
    05C50, 
    05C75, 
    31C20, 
    46E22, 
    47B25, 
    47B32, 
    60J10, 
    Secondary:
    31C35, 
    47B39, 
    82C41. 
    }


  \date{\bf\today. \q \version{}{Rough version (for editing).}}

\maketitle

\setcounter{tocdepth}{1}
{\small \tableofcontents}

\allowdisplaybreaks

\input{introduction}

\input{setup}

\input{discrete-gauss-green-formula}

\input{monopoles}

\input{resistance-metrics}

\input{gelfand}

\input{wiener}

\input{boundary}

\input{examples}

\subsection*{Acknowledgements}
The authors are grateful to Ecaterina Sava and Wolfgang Woess for organizing the Boundaries 2009 workshop and to Florian Sobieczky for organizing the 2009 Alp Workshop on Spectral Theory and Random Walks, thereby making it possible for us to meet and exchange ideas with a multitude of top-tier researchers. We are grateful to the  participants of these workshops for their ideas and comments, suggestions and general mathematical stimulation (there were a lot of great talks!). In particular, we benefitted from conversations with Donald Cartwright, Massimo Picardello, Agelos Georgakopolous, Vadim Kaimanovich, Matthias Keller, Elmar Teufl, Wolfgang Woess, and Radek Wojciechowski.

\bibliographystyle{alpha}
\bibliography{RBIN}

\end{document}

%% file: introduction.tex

\section{Introduction}
\label{sec:introduction}

There are several notions of ``boundary'' as ``points at infinity'' associated to infinite graphs. Some of these come directly from graph theory, like the notion of \emph{graph ends} \cite{PicWoess90,Diestel2006} or \emph{ideal boundary} \cite{MuYa89,Yamasaki86}. Others come from by way of the associated reversible Markov process, the random walk associated to the graph, like topological notion of \emph{Martin boundary} \cite{PicWoess87,Saw97} or its measure-theoretic refinement, the Poisson boundary \cite{KarlssonWoess,Kaimanovich91}. There are also less well-known ideas, like the discrete Royden boundary \cite{KayanoYamasaki88} and discrete Kuramochi boundary \cite{MurakamiYamasaki}. Interrelations amongst these concepts are detailed in two excellent collections of notes: \cite{Woess00} and \cite{Soardi94}. This material has its roots in minimal surface theory, probability theory, ergodic theory, and group theory, and the central ideas are often analogues of a corresponding notion for continuous domains (manifolds, Lie groups, etc.). 

This paper gives a brief account of a new type of boundary developed in \cite{DGG,ERM,bdG} and \cite{OTERN} which we call the \emph{resistance boundary}; it is denoted $\bd G$. It bears many similarities to the Martin and Poisson boundaries, but pertains to a different class of functions: the functions of \emph{finite energy}.
Let $G$ be a \emph{resistance network} (i.e., a connected simple weighted graph) with vertex set \verts and edges determined by a symmetric \emph{conductance} function \cond which weights the edges: $\cond_{xy} = \cond_{yx} \geq 0$, and $\cond_{xy} > 0$ iff there is an edge from $x$ to $y$, which is denoted $x \nbr y$. The energy of a function $u:G \to \bC$ is then defined to be 
\linenopax
\begin{align}\label{eqn:intro-energy}
  \energy(u) :=  \frac12 \sum_{x,y \in \verts} \cond_{xy} |u(x)-u(y)|^2.
\end{align}
For the most part, it suffices to work with \bR-valued functions (see Remark~\ref{rem:Real-and-complex-valued-functions-on-verts} in particular). However, we will need \bC-valued functions for some applications of spectral theory in \S\ref{sec:boundary-integral-representation-for-harmonics}.

Under suitable hypotheses, if $h$ is a bounded harmonic function on $X$, then Poisson boundary theory provides a measure space $(\del X, \gm)$ with respect to which one has an integral representation of $h$ in terms of a kernel $k:X\times \del X \to \bC$:
\linenopax
\begin{align}\label{eqn:Poisson-rep-preview}
  h(x) = \int_{\del X} k(x,\gx) \tilde{h}(\gx) \,d\gm(\gx),
\end{align}
where $\tilde{h}$ is the extension of $h$ to $\del X$, in some sense. This paper provides a synopsis of how one can obtain a similar representation for the harmonic functions of finite energy. However, instead of using ergodic theory or (topological) compactifications, we take an entirely different approach: operator theory and functional analysis. 

After embedding the resistance network into a certain Hilbert space, we construct a space of distributions (i.e., generalized functions) on that Hilbert space. We then show that this space of distributions contains the boundary of the original network, in the sense that it supports integral representations of harmonic functions on the network. We work with the \emph{energy space}, a Hilbert space whose inner product is given by the sesquilinear form associated to \energy by polarizing \eqref{eqn:intro-energy}:
\linenopax
\begin{align}\label{eqn:intro-energy-form}
  \la u,v \ra_\energy 
  = \energy(u,v) 
  := \frac12 \sum_{x,y \in \verts} \cond_{xy} \left(\cj{u(x)}-\cj{u(y)}\right)\left(v(x) - \vstr[2.2]v(y)\right).
\end{align}
We construct a reproducing kernel for this Hilbert space, and then use it to obtain a Gel'fand triple
\linenopax
\begin{align}\label{eqn:Gelfand-triple-preview}
  \Schw \ci \HE \ci \Schw'.
\end{align}
Here, \Schw is a \energy-dense subspace of \HE which is also equipped with a strictly finer ``test function topology'' (defined in terms of the domain of the Laplacian), and
the space $\Schw'$ is the dual space of \Schw with respect to this finer topology; the specifics are discussed further just below. For now, however, let us eschew technical details and just say that $\Schw'$ is strictly larger than \HE, and it is in $\Schw'$ that the boundary $\bd G$ lies. This framework allows us to invoke Minlos' theorem and Wiener's isometric embedding theorem, powerful tools from the theory of stochastic integration. Boundary theory usually involves an enlargement of the original space, either by topological means (e.g., by compactification or completion, in the case of Martin boundary) or by measure-theoretical means (e.g., by taking the measurable hull of an equivalence relation, as in Poisson boundary). For the resistance boundary $\bd \Graph$, we enlarge \HE (the Hilbert space representation of the resistance network) by embedding it into $\Schw'$ via the inclusion map.


\begin{defn}\label{def:graph-laplacian}
  The \emph{Laplacian} on a resistance network $(\Graph,\cond)$ is the linear difference operator \Lap which acts on a function $v:\verts \to \bC$ by
  \linenopax
  \begin{equation}\label{eqn:def:laplacian}
    (\Lap v)(x) :
    = \sum_{y \nbr x} \cond_{xy}(v(x)-v(y)).
  \end{equation}
  A function $v:\verts \to \bC$ is \emph{harmonic} iff $\Lap v(x)=0$ for each $x \in \verts$. 
\end{defn}

Note that we adopt the (physicists') sign convention in \eqref{eqn:def:laplacian} (so that the spectrum is nonnegative) and thus our Laplacian is the negative of the one commonly found in the PDE literature; e.g., \cite{Kig01,Str06}. 

The study of resistance boundaries begins with the following well-known identity for finite networks.
\begin{prop}\label{thm:E(u,v)=<u,Lapv>}
  Let \Graph be a finite network. For functions $u,v$ on \verts,
  \linenopax
  \begin{equation}\label{eqn:E(u,v)=<u,Lapv>}
    \energy(u,v)
    = \sum_{x \in \verts} \cj{u}(x) \Lap v(x).
  \end{equation}
\end{prop}
The right-hand side of \eqref{eqn:E(u,v)=<u,Lapv>} is often denoted by $\la u, \Lap v\ra_2$. Theorem~\ref{thm:E(u,v)=<u,Lapv>+sum(normals)} gives a broad extension of Proposition~\ref{thm:E(u,v)=<u,Lapv>} to a certain domain \MP (see Definition~\ref{def:LapM}). 
Extensions of this type have been studied before (see \cite{Maeda80,Kayano88}), but only with regard to determining conditions that ensure $\energy(u,v) = \la u, \Lap v\ra_2$. By contrast, we are more interested in the situation for which it is replaced by 
\linenopax
  \begin{equation}\label{eqn:DGG-preview}
    \la u, v \ra_\energy
    = \sum_{\verts} \cj{u} \Lap v
      + \sum_{\bd \Graph} \cj{u} \dn v.
  \end{equation}
Theorem~\ref{thm:E(u,v)=<u,Lapv>+sum(normals)} gives conditions under which \eqref{eqn:DGG-preview} holds; the notation $\bd \Graph$ and $\dn v$ are explained precisely in Definition~\ref{def:subgraph-boundary} and Definition~\ref{def:boundary-sum}. In particular, \eqref{eqn:DGG-preview} holds for any $u \in \HE$ when $v$ lies in a certain dense subspace of \HE which we denote by \MP. The space \MP was introduced in \cite{DGG} for this purpose and also to serve as a dense domain for the possibly unbounded Laplace operator, which will be useful later for the construction of \Schw. 
We call \eqref{eqn:DGG-preview} the discrete Gauss-Green identity by analogy with
\linenopax
  \begin{equation*}
    \int_{\gW} \grad u \grad v \,dV
    = -\int_{\gW} u \Lap v \,dV
      + \int_{\del\gW} u \dn v \,dS.
  \end{equation*}
 
The space \HE consists of potentials (functions on the vertices of \Graph, modulo constants; see Definition~\ref{def:The-energy-Hilbert-space}) and enjoys an orthogonal decomposition into the subspace \Fin of finitely supported functions and the subspace \Harm of harmonic functions; this is given precisely in Definitions~\ref{def:Fin}--\ref{def:Harm} and Theorem~\ref{thm:HE=Fin+Harm}. It turns out that \HE has a reproducing kernel $\{v_x\}_{x \in \verts}$: for any $u \in \HE$, one has 
\linenopax
\begin{align*}
  \la v_x, u \ra_\energy = u(x) - u(o), \qq \forall x \in \verts,
\end{align*}
where $o \in \verts$ is a fixed reference point. Since the reproducing kernel behaves well with respect to (orthogonal) projections $P$, we also have reproducing kernels $\{f_x\}_{x \in \verts}$ for \Fin and $\{h_x\}_{x \in \verts}$ for \Harm, where
\linenopax
\begin{align*}
  f_x := \Pfin v_x,
  \q\text{and}\q
  h_x := \Phar v_x.
\end{align*}

In Theorem~\ref{thm:Boundary-representation-of-harmonic-functions}, we apply \eqref{eqn:DGG-preview} to the reproducing kernels $\{h_x\}_{x \in \verts}$ for \Harm, and find that for all $h \in \Harm$, 
  \linenopax
  \begin{equation}\label{eqn:boundary-repn-for-harmonic}
    h(x) - h(o) = \sum_{\bd \Graph} h \dn{h_x}.
  \end{equation}
This direct analogue of \eqref{eqn:Poisson-rep-preview} first appeared in \cite[Cor.~3.14]{DGG}. 
Formula \eqref{eqn:boundary-repn-for-harmonic} gives a boundary sum representation of harmonic functions, but the boundary sum in \eqref{eqn:boundary-repn-for-harmonic} is understood only as a limit of sums taken over boundaries of finite subnetworks. Comparison of \eqref{eqn:boundary-repn-for-harmonic} and \eqref{eqn:Poisson-rep-preview} makes one optimistic that $\bd \Graph$ can be realized as a measure space which supports a measure corresponding to $\dn{h_x}$, thus replacing the sum in \eqref{eqn:boundary-repn-for-harmonic} with a integral. In Corollary~\ref{thm:Boundary-integral-repn-for-harm}, we extend \eqref{eqn:boundary-repn-for-harmonic} to such an integral representation for which \eqref{eqn:boundary-repn-for-harmonic} is analogous to a Riemann sum.

The primary difference between our boundary theory and that of Poisson and Martin is rooted in our focus on \HE: both of these classical theories concern harmonic functions with growth/decay restrictions. By contrast, provided they neither grow too wildly nor oscillate too wildly, elements of \HE may be unbounded and may fail to remain nonnegative. 
From \cite{AnconaLyonsPeres}, it is known that functions which are \energy-limits of finitely supported functions must vanish at \iy (except for a set of measure 0 with respect to the usual path-space measure); however see \cite[Ex.~13.10]{OTERN} for an unbounded harmonic function of finite energy. Note, however, that functions of finite energy can always be approximated in \HE 
by bounded functions; cf.~\cite[\S3.7]{Soardi94}.

Just as for Martin and Poisson boundaries, the resistance boundary essentially consists of different limiting behaviors of the (transient) random walk on the network, as the walker tends to infinity. 
It turns out that 
recurrent networks have no resistance boundary, and transient networks with no nontrivial harmonic functions have exactly one boundary point (corresponding to the fact that the monopole at $x$ is unique; see Definition~\ref{def:monopole}). In particular, the integer lattices $(\bZd,\one)$ each have 1 boundary point for $d \geq 3$ and 0 boundary points for $d=1,2$. Further examples are discussed in \S\ref{sec:examples}. 


\subsubsection*{Outline}

\S\ref{sec:electrical-resistance-networks} recalls basic definitions and some previously obtained results. In particular, we give precise definitions for the Laplace operator \Lap, the energy space \HE, the reproducing kernel $\{v_x\}$, monopoles $w_x$, the monopolar domain \MP, and we discuss the Royden decomposition of \HE into the finitely supported functions and the harmonic functions.
\S\ref{sec:relating-energy-form-to-Laplacian} states the discrete Gauss-Green identity and gives the definition of the boundary sum $\sum_{\bd G} u \dn v$, as a limit of sums. Some implications of the discrete Gauss-Green identity are given, including several characterizations of transience of the random walk on the network.
\S\ref{sec:Effective-resistance} gives the definition of effective resistance, and discusses how this metric can be extended to infinite networks in different ways 
the \emph{free resistance} $R^F(x,y)$ and \emph{wired resistance} $R^W(x,y)$. 
\S\ref{sec:boundary-integral-representation-for-harmonics} discusses the boundary sum representation for elements of \Harm as introduced in \eqref{eqn:boundary-repn-for-harmonic}. This section also gives an overview of the theory of Gel'fand triples, Minlos' theorem, and Wiener's theorem, and how these enable one to obtain a Gaussian probability measure on the space $\Schw'$ alluded to in \eqref{eqn:Gelfand-triple-preview}.
\S\ref{sec:boundary-integral-representation-for-harmonics} 
gives the boundary integral representation of elements of \Harm: an integral version of \eqref{eqn:boundary-repn-for-harmonic} which is an \HE-analogue of \eqref{eqn:Poisson-rep-preview}. 
\S\ref{sec:examples} contains several examples which illustrate our results.

\pgap

Boundary theory is a well-established subject; the deep connections between harmonic analysis, probability, and potential theory have led to several notions of boundary and we will not attempt to give complete references. However, we recommend \cite{Saw97} for introductory material on Martin boundary and \cite{Woess00} for a more detailed discussion. Introductory material on resistance networks may be found in \cite{DoSn84} and \cite{Lyons:ProbOnTrees}, and \cite{Kig03} gives a detailed investigation of resistance forms (a potential-theoretic generalization of resistance networks). 
More specific background appears in \cite{TerryLyons, Car73a} and the foundational paper \cite{Nash-Will59}. With regard to infinite graphs and finite-energy functions, see \cite{Soardi94, Woess00, SoardiWoess91, CaW92, Dod06, PicWoess90, PicWoess88, Wo86, Thomassen90}. Applications to analysis on fractals can be found in \cite{Kig01, Str06}. For papers studying fractals as boundaries of networks or Markov processes, see \cite{DenkerSato99, DenkerSato01, DenkerSato02, Kaimanovich01, LauWang09, JuLauWang10, Kig09}.

%% file: setup.tex

\section{The energy space \HE} 
\label{sec:The-energy-space}
\label{sec:electrical-resistance-networks}

We now proceed to introduce the key notions used throughout this paper: resistance networks, the energy form \energy, the Laplace operator \Lap, the energy space \HE, the reproducing kernel $\{v_x\}$, and their elementary properties. 

\begin{defn}\label{def:ERN}
  A resistance network is a connected graph $(\Graph,\cond)$, where \Graph is a graph with vertex set \verts, and \cond is the \emph{conductance function} which defines adjacency by $x \nbr y$ iff $c_{xy}>0$, for $x,y \in \verts$. We assume $\cond_{xy} = \cond_{yx} \in [0,\iy)$, and write $\cond(x) := \sum_{y \nbr x} \cond_{xy}$. We require $\cond(x) < \iy$ but $\cond(x)$ need not be a bounded function on \verts, and note that vertices of infinite degree are allowed. The notation \cond may be used to indicate the multiplication operator $(\cond v)(x) := \cond(x) v(x)$, 
  \version{}{\marginpar{Is it ok to say ``basis'' here?}}
  i.e., the diagonal matrix with entries $\cond(x)$ with respect to the (vector space) basis $\{\gd_x\}$.
\end{defn}

As the letters $x,y,z$ always denote vertices, it causes no confusion to write $x,y,z \in \Graph$ instead of $x,y,z \in \verts$. Similarly, $u$ and $v$ will always denote functions which map vertices to scalars, e.g., $u:\verts \to \bC$.

In Definition~\ref{def:ERN}, ``connected'' means simply that for any $x,y \in \Graph $, there is a finite sequence $\{x_i\}_{i=0}^n$ with $x=x_0$, $y=x_n$, and $\cond_{x_{i-1} x_i} > 0$, $i=1,\dots,n$. Conductance is the reciprocal of resistance, so one can think of $(\Graph,\cond)$ as a network of nodes \verts connected by resistors of resistance $\cond_{xy}^{-1}$. We may assume there is at most one edge from $x$ to $y$, as two conductors $\cond^1_{xy}$ and $\cond^2_{xy}$ connected in parallel can be replaced by a single conductor with conductance $\cond_{xy} = \cond^1_{xy} + \cond^2_{xy}$. Also, we assume $\cond_{xx}=0$ so that no vertex has a loop, as electric current will never flow along a conductor connecting a node to itself. 

\begin{defn}\label{def:exhaustion-of-G}
  An \emph{exhaustion} of \Graph is an increasing sequence of finite and connected subgraphs $\{\Graph_k\}_{k=1}^\iy$, so that $\Graph_k \ci \Graph_{k+1}$ and $\Graph = \bigcup \Graph_k$. Since any vertex or edge is eventually contained in some $G_k$, there is no loss of generality in assuming they are contained in $G_1$, for the purposes of a specific computation.
\end{defn}

\begin{defn}\label{def:infinite-vertex-sum}  
  The notation
  \linenopax
  \begin{equation}\label{eqn:def:infinite-sum}
    \sum_{x \in \Graph} := \lim_{k \to \iy} \sum_{x \in \Graph_k}
  \end{equation}
  is used whenever the limit is independent of the choice of exhaustion $\{\Graph_k\}$ of \Graph. This is clearly justified, for example, whenever the sum has only finitely many nonzero terms, or is absolutely convergent as in the definition of \energy just below.
\end{defn}

\begin{defn}\label{def:graph-energy}
  The \emph{energy} of functions $u,v:\verts \to \bC$ is given by the (closed, bilinear) Dirichlet form
  \linenopax
  \begin{align}\label{eqn:def:energy-form}
    \energy(u,v)
    := \frac12 \sum_{x \in \Graph}  \sum_{y \in \Graph} \cond_{xy}(\cj{u}(x)-\cj{u}(y))(v(x)-v(y)),
  \end{align}
  with the energy of $u$ given by $\energy(u) := \energy(u,u)$.
  The \emph{domain of the energy} is
  \linenopax
  \begin{equation}\label{eqn:def:energy-domain}
    \dom \energy = \{u:\verts \to \bC \suth \energy(u)<\iy\}.
  \end{equation}
\end{defn}

Since $\cond_{xy}=\cond_{yx}$ and $\cond_{xy}=0$ for nonadjacent vertices, the initial factor of $\frac12$ in \eqref{eqn:def:energy-form} implies there is exactly one term in the sum for each edge in the network.

\begin{defn}\label{def:H_energy}\label{def:The-energy-Hilbert-space}
  Let \one denote the constant function with value 1 and recall that $\ker \energy = \bC \one$. 
 The energy form \energy is symmetric and positive definite on $\dom \energy$. Then $\dom \energy / \bC \one$ is a vector space with inner product and corresponding norm given by
  \linenopax
  \begin{equation}\label{eqn:energy-inner-product}
    \la u, v \ra_\energy := \energy(u,v)
    \q\text{and}\q
    \|u\|_\energy := \energy(u,u)^{1/2}.
  \end{equation}
  The \emph{energy Hilbert space} is defined to be 
  \linenopax
  \begin{align}\label{eqn:HE}
    \HE := \frac{\dom \energy}{\ker \energy} = \frac{\dom \energy}{\bC \one}.
  \end{align}
\end{defn}

Thus, \HE consists of potentials: we are not interested in values $u(x)$ as much as differences $u(x)-u(y)$. In other words, if $u$ and $v$ are both elements of $\dom \energy$ and there is some constant $k \in \bC$ such that $u(x)-v(x) = k$ for all $x \in \Graph$, then $u$ and $v$ are both representatives of the same element (equivalence class) of \HE.

\begin{defn}\label{def:vx}\label{def:energy-kernel}
  Let $v_x$ be defined to be the unique element of \HE for which
  \linenopax
  \begin{equation}\label{eqn:def:vx}
    \la v_x, u\ra_\energy = u(x)-u(o),
    \qq \text{for every } u \in \HE.
  \end{equation}
  The collection $\{v_x\}_{x \in \Graph}$ forms a reproducing kernel for \HE; cf.~\cite[Cor.~2.7]{DGG}. We call it the \emph{energy kernel} and \eqref{eqn:def:vx} shows its span is dense in \HE. Note that $v_o$ corresponds to a constant function, since $\la v_o, u\ra_\energy = 0$ for every $u \in \HE$. Therefore, $v_o$ (or $o$) may often ignored or omitted.
\end{defn} 

\begin{defn}\label{def:dipole}
  A \emph{dipole} is any $v \in \HE$ satisfying the pointwise identity $\Lap v = \gd_x - \gd_y$ for some vertices $x,y \in \Graph$. 
  The elements of the energy kernel are all dipoles: one can check that $\Lap v_x = \gd_x - \gd_o$ as in \cite[Lemma~2.13]{DGG}.
\end{defn}

\begin{remark}\label{rem:G-as-index}
  To minimize cumbersome notation, let $\{x \in \Graph\}$ be the default index set from now on. That is, we use $\{v_x\}$ to denote the energy kernel $\{v_x\}_{x \in \Graph}$, and $\spn\{v_x\}$ to denote the set of all linear combinations of elements of $\{v_x\}$, etc.
\end{remark}

\subsection{The finitely-supported functions and the harmonic functions}

\begin{defn}\label{def:Fin}
  For $v \in \HE$, one says that $v$ has \emph{finite support} iff there is a finite set $F \ci \verts$ for which $v(x) = k \in \bC$ for all $x \notin F$, i.e., the set of functions of finite support in \HE is $\spn\{\gd_x\}$, 
  where $\gd_x$ is the Dirac mass at $x$, i.e., the element of \HE containing the characteristic function of the singleton $\{x\}$. 
  Define \Fin to be the closure of $\spn\{\gd_x\}$ with respect to \energy. 
\end{defn}

\begin{remark}\label{rem:Diracs-not-onb}
  The usual candidate for an orthonormal basis (onb) in $\ell^2(\verts)$ would be the collection of Dirac masses $\{\gd_x\}$; however, this is not an onb in \HE. One can compute from \eqref{eqn:def:energy-form} that 
  \linenopax
  \begin{align}\label{eqn:E(Dirac)=cond}
     \la \gd_x, \gd_y\ra_\energy = \energy(\gd_x,\gd_y) = -\cond_{xy}, \q\text{ for } x \neq y,
     \qq\text{and}\qq
     \energy(\gd_x) = \cond(x),
  \end{align}
  so that $\gd_x \not\perp \gd_y$ with respect to \energy.
  Moreover, Theorem~\ref{thm:HE=Fin+Harm} shows that $\{\gd_x\}$ is, in general, \emph{not even dense in \HE.} It is immediate from \eqref{eqn:E(Dirac)=cond} that $\gd_x \in \HE$.
\end{remark}

\begin{defn}\label{def:Harm}
  The harmonic subspace of \HE is denoted
  \linenopax
  \begin{equation}\label{eqn:Harm}
    \Harm := \{v \in \HE \suth \Lap v(x) = 0, \text{ for all } x \in \Graph\}.
  \end{equation}
  Note that this is independent of choice of representative for $v$ in virtue of \eqref{eqn:def:laplacian}.
\end{defn}

The following result is sometimes called the ``Royden Decomposition'' since \cite[Thm.~4.1]{Yamasaki79}, in reference to Royden's analogous result for Riemann surfaces; see \cite[\S{VI}]{Soardi94}, \cite[\S9.3]{Lyons:ProbOnTrees}\footnote{This name is also sometimes associated with the corresponding (nonorthogonal) decomposition for the ``grounded energy form''; see Remark~\ref{rem:grounded-energy}.}. It follows immediately from \cite[Lemma~2.11]{DGG}, which states that $\la \gd_x, u \ra_\energy = \Lap u(x)$ for any $x \in \Graph$; cf.~\cite[Thm.~2.15]{DGG}.

\begin{theorem}[Royden decomposition]
  \label{thm:HE=Fin+Harm}
  $\HE = \Fin \oplus \Harm$.
\end{theorem}

\begin{defn}\label{def:ux}
  Let $f_x = \Pfin v_x$ denote the image of $v_x$ under the (orthogonal) projection to \Fin. Similarly, let $h_x = \Phar v_x$ denote the image of $v_x$ under the projection to \Harm. 
\end{defn}

\begin{remark}[Reproducing kernels for \Fin and \Harm]
  \label{rem:Real-and-complex-valued-functions-on-verts}
The reproducing kernel property behaves well with respect to orthogonal projections, and consequently, $\{f_x\}$ is a reproducing kernel for \Fin, and $\{h_x\}$ is a reproducing kernel for \Harm. 
  While we will need complex-valued functions for some results obtained via spectral theory, it will usually suffice to consider \bR-valued functions because the reproducing kernels elements $v_x, f_x, h_x$ all have \bR-valued representatives \cite[Lemma~2.24]{DGG}.
\end{remark}

\subsection{Monopoles}

\begin{defn}\label{def:monopole}
  A \emph{monopole} at $x \in \Graph$ is an element $w_x \in \HE$ which satisfies
  $\Lap w_x(y) = \gd_{xy}$, where $\gd_{xy}$ is Kronecker's delta. 
  In case the network supports monopoles (that is, if the above Dirichlet equation admits finite-energy solutions), let $w_o$ always denote the unique energy-minimizing monopole at the origin. 
  
  With $v_x$ and $f_x = \Pfin v_x$ as above, we indicate the distinguished monopoles 
  \linenopax
  \begin{align}\label{eqn:def:monov-and-monof}
    \monov := v_x + w_o
    \q\text{and}\q
    \monof := f_x + w_o.
  \end{align}
\end{defn}

\begin{remark}\label{rem:w_O-in-Fin}
  Note that $w_o \in \Fin$, whenever it is present in \HE, and similarly that \monof is the energy-minimizing monopole at $x$. To see this, suppose $w_x$ is any monopole at $x$. Since $w_x \in \HE$, write $w_x = f+h$ by Theorem~\ref{thm:HE=Fin+Harm}, and get $\energy(w_x) = \energy(f) + \energy(h)$. Projecting away the harmonic component will not affect the monopole property, so $\monof = \Pfin w_x$ is the unique monopole of minimal energy.  The Green function is $g(x,y) = w_y^o(x)$, where $w_y^o$ is the representative of \monof[y] which vanishes at \iy. 
\end{remark}

\begin{defn}\label{def:LapM}
  The dense subspace of \HE spanned by monopoles (and dipoles) is
  \linenopax
  \begin{equation}\label{eqn:def:MP}
    \MP := \spn\{v_x\} + \spn\{\monov, \monof\}.
  \end{equation}
  Let \LapM be the closure of the Laplacian when taken to have the dense domain \MP. \end{defn}

Since \Lap agrees with \LapM pointwise, we may suppress reference to the domain for ease of notation. It is shown in \cite[Lemma~3.5]{DGG} that \LapM is Hermitian with $\la u, \LapM u\ra_\energy \geq 0$ for all $u \in \MP$. When given a pointwise identity $\Lap u = v$, there is an associated identity in \HE, but one must use the adjoint: $\Lap u(x) = v(x)$ for all $x \in \Graph$ if and only if $v = \LapM^\ad u$ in \HE \cite[Lemma~3.7]{DGG}. Note that \LapM may have \emph{defect vectors}; such an object is an element of the Hilbert space \HE (though clearly not an element of \MP) which has a representative $u$ satisfying
\linenopax
\begin{align*}
  \LapM u(x) = -u(x), \forall x \in \Graph,
  \qq \text{and} \qq
  u \in \dom \LapM^\ad.
\end{align*}
See \cite[\S4.2]{SRAMO} or \cite[\S13.4]{OTERN}. While it is always the case that the (possibly unbounded) operator \LapM is Hermitian (i.e. $\LapM \ci \LapM^\ad$), this shows that \LapM may fail to be self-adjoint (i.e. $\LapM = \LapM^\ad$).

\begin{remark}[Monopoles and transience]
  \label{rem:transient-iff-monopoles}
  The presence of monopoles in \HE is equivalent to the transience of the simple random walk on the network with transition probabilities $p(x,y) = \cond_{xy}/\cond(x)$: note that if $w_x$ is a monopole, then the current induced by $w_x$ is a unit flow to infinity with finite energy. It was proved in \cite{TerryLyons} that the network is transient if and only if there exists a unit current flow to infinity; see also \cite[Thm.~2.10]{Lyons:ProbOnTrees}. 
  Moreover, it is shown in \cite[Lemma~3.6]{DGG} that when the network is transient, \MP contains the spaces $\spn\{v_x\}, \spn\{f_x\}$, and $\spn\{h_x\}$, where $f_x = \Pfin v_x$ and $h_x = \Phar v_x$. When $\Harm = 0$ (in particular, when the network is not transient), $f_x = v_x$ and so $\MP = \spn\{v_x\}= \spn\{f_x\}$ trivially.
\end{remark}

%% file: discrete-gauss-green-formula.tex

\section{The discrete Gauss-Green formula}
\label{sec:relating-energy-form-to-Laplacian}

In Theorem~\ref{thm:E(u,v)=<u,Lapv>+sum(normals)}, we establish a discrete version of the Gauss-Green formula which extends Proposition~\ref{thm:E(u,v)=<u,Lapv>} to the case of infinite graphs; the scope of validity of this formula is given in terms of the space \MP of Definition~\ref{def:monopole}. The appearance of a somewhat mysterious boundary term alluded to in \eqref{eqn:DGG-preview} prompts several questions which are discussed in Remark~\ref{rem:boundary-term}. 


\subsection{Relating \Lap to \energy}
\label{sec:Relating-Lap-to-energy}

\begin{defn}\label{def:subgraph-boundary}
  If $H$ is a subgraph of $G$, then the boundary of $H$ is
  \linenopax
  \begin{equation}\label{eqn:subgraph-boundary}
    \bd H := \{x \in H \suth \exists y \in H^\complm, y \nbr x\}.
  \end{equation}
  The \emph{interior} of a subgraph $H$ consists of the vertices in $H$ whose neighbours also lie in $H$:
  \linenopax
  \begin{equation}\label{eqn:interior}
    \inn H := \{x \in H \suth y \nbr x \implies y \in H\} = H \less \bd H.
  \end{equation}
  For vertices in the boundary of a subgraph, the \emph{normal derivative} of $v$ is
  \linenopax
  \begin{equation}\label{eqn:sum-of-normal-derivs}
    \dn v(x) := \sum_{y \in H} \cond_{xy} (v(x) - v(y)),
    \qq \text{for } x \in \bd H.
  \end{equation}
  Thus, the normal derivative of $v$ is computed like $\Lap v(x)$, except that the sum extends only over the neighbours of $x$ which lie in $H$.
\end{defn}
  \glossary{name={$\dn v$},description={normal derivative of a function with respect to a subgraph},sort=d,format=textbf}

Definition~\ref{def:subgraph-boundary} will be used primarily for subgraphs that form an exhaustion of \Graph, in the sense of Definition~\ref{def:exhaustion-of-G}. 

\begin{defn}\label{def:boundary-sum}
  A \emph{boundary sum} is computed in terms of an exhaustion $\{G_k\}$ by
  \linenopax
  \begin{equation}\label{eqn:boundary-sum}
    \sum_{\bd \Graph} := \lim_{k \to \iy} \sum_{\bd \Graph_k},
  \end{equation}
  whenever the limit is independent of the choice of exhaustion, as in Definition~\ref{def:infinite-vertex-sum}.
\end{defn}

\begin{theorem}[Discrete Gauss-Green Formula]
  \label{thm:E(u,v)=<u,Lapv>+sum(normals)}
  If $u \in \HE$ and $v \in \MP$, then
  \linenopax
  \begin{equation}\label{eqn:E(u,v)=<u_0,Lapv>+sum(normals)}
    \la u, v \ra_\energy
    = \sum_{\Graph} \cj{u} \Lap v
      + \sum_{\bd \Graph} \cj{u} \dn v.
  \end{equation}
\end{theorem}

\begin{cor}\label{thm:sumLap(u)=-sumdn(u)}
  For all $u \in \dom \LapM$, $\sum_{\Graph} \Lap u = - \sum_{\bd \Graph} \dn u$. Thus, the discrete Gauss-Green formula \eqref{eqn:E(u,v)=<u_0,Lapv>+sum(normals)} is independent of choice of representatives.
\end{cor}

\begin{remark}\label{rem:DGG-holds-without-hypotheses}
  The proof of Theorem~\ref{thm:E(u,v)=<u,Lapv>+sum(normals)} follows from taking limits of 
  \linenopax
  \begin{align*}
    \sum_{x \in G_k} \cj{u}(x) \Lap v(x)
      + \sum_{x \in \bd G_k} \cj{u}(x) \dn v(x).
  \end{align*}
  Thus, the decomposition \eqref{eqn:E(u,v)=<u_0,Lapv>+sum(normals)} is true for all $u,v \in \HE$, but is meaningless if it takes the form $\iy - \iy$. A key point of Theorem~\ref{thm:E(u,v)=<u,Lapv>+sum(normals)} is that for $u,v$ in the specified domains, the two sums are both finite and independent of choice of exhaustion. However, the specific value of each sum \emph{is} dependent on the choice of representative for $u$; this motivates Definition~\ref{def:boundary-term-is-nonvanishing}.
  
  It is also clear that \eqref{eqn:E(u,v)=<u_0,Lapv>+sum(normals)} remains true much more generally than under the specified conditions; certainly the formula holds whenever $\sum_{x \in \Graph} \left|{u}(x) \Lap v(x)\right| < \iy$.
  Unfortunately, given any hypotheses more specific than this, the limitless variety of infinite networks almost always allows one to construct a counterexample; i.e. one cannot give a condition for which the formula is true for all $u \in \HE$, for all networks. 
  Nonetheless, the formula remains true and even useful in many specific and general contexts. For example, it is clearly valid whenever $v$ is a dipole, including all those in the energy kernel. We will also see that it holds for the projections of $v_x$ to \Fin and to \Harm. Consequently, for $v$ which are limits of elements in \MP, we can use this result in combination with ad hoc arguments. 
  
  A formula similar to \eqref{eqn:E(u,v)=<u_0,Lapv>+sum(normals)} appears in \cite[Prop~1.3]{DodziukKarp88}; however, these authors apparently do not pursue the extension of this formula to infinite networks. Another similar result appears in \cite[Thm.~4.1]{Kayano88}, where the authors give some conditions under which \eqref{eqn:E(u,v)=<u,Lapv>} extends to infinite networks. The main differences here are that the scope of Kayano and Yamasaki's theorem is limited to a subset of what we call \Fin, and that Kayano and Yamasaki are interested in when the boundary term vanishes; we are more interested in when it is finite and nonvanishing; see Theorem~\ref{thm:transience}, for example. Since Kayano and Yamasaki do not discuss the structure of the space of functions they consider, it is not clear how large the scope of their result is; their result requires the hypothesis $\sum_{x \in \Graph} \left|{u}(x) \Lap v(x)\right| < \iy$, but it is not so clear what functions satisfy this. By contrast, we develop a dense subspace of functions on which to apply the formula. Furthermore, in the forthcoming paper \cite{ERM}, we show that these functions are relatively easy to compute.
\end{remark}

\begin{remark}\label{rem:boundary-term}
  We refer to $\sum_{\bd \Graph} u \dn v$ as the ``boundary term'' by analogy with classical PDE theory. This terminology should not be confused with the notion of boundary that arises in the discussion of the discrete Dirichlet problem, where the boundary is a prescribed subset of \verts.   As the boundary term may be difficult to contend with, it is extremely useful to know when it vanishes, for example:
  \begin{enumerate}[(i)]
    \item when the network is recurrent (Theorem~\ref{thm:transience}),
    \item when $v$ is an element of the energy kernel \cite[Lemma~5.8]{DGG},
    \item when $u,v,\Lap u, \Lap v$ lie in $\ell^2$ \cite[Lemma~5.12]{DGG}, and
    \item when either $u$ or $v$ has finite support \cite[Lemma~5.16]{DGG}.
  \end{enumerate}
\end{remark}

%% file: monopoles.tex

\subsection{More about monopoles and the space \MP}
\label{sec:More-about-monopoles}

This section studies the role of the monopoles with regard to the boundary term of Theorem~\ref{thm:E(u,v)=<u,Lapv>+sum(normals)}, and provides several characterizations of transience of the network, in terms the operator-theoretic properties of \LapM. 
  
Note that if $h \in \Harm$ satisfies the hypotheses of Theorem~\ref{thm:E(u,v)=<u,Lapv>+sum(normals)}, then $\energy(h) = \sum_{\bd\Graph} h \dn h$. On the other hand, $\energy(u) = \sum_{\Graph} u \Lap u$ for all $u \in \HE$ iff the network is recurrent, as stated in Theorem~\ref{thm:transience}. With respect to $\HE = \Fin \oplus \Harm$, this shows that the energy of finitely supported functions comes from the sum over \Graph, and the energy of harmonic functions comes from the boundary sum. However, for a monopole $w_x$, the representative specified by $w_x(x)=0$ satisfies $\energy(w) = \sum_{\bd\Graph} w \dn w$ but the representative specified by $w_x(x) = \energy(w_x)$ satisfies $\energy(w) = \sum_{\Graph} w \Lap w$. Roughly speaking, a monopole is therefore ``half of a harmonic function'' or halfway to being a harmonic function. A further justification for this comment is given by Corollary~\ref{thm:Harm-nonzero-iff-multiple-monopoles} (the proof shows that a harmonic function can be constructed from two monopoles at the same vertex, see \cite[Cor.~4.4]{DGG}). The general theme of this section is the ability of monopoles to ``bridge'' the finite and the harmonic.

\begin{theorem}[{\cite[Thm.~1.33]{Soardi94}}]
  \label{thm:Soardi's-harmonic-implies-transient}
  Let $u$ be a nonnegative function on a recurrent network. Then $u$ is superharmonic if and only if $u$ is constant.  
\end{theorem}

It follows from Theorem~\ref{thm:Soardi's-harmonic-implies-transient} that $\Harm \neq 0$ implies the existence of a monopole in \HE, i.e., the transience of the network; cf.~ \cite[Cor.~4.3]{DGG}. However, it turns out that a nontrivial harmonic function can only exist when there is more than one monopole.

\begin{cor}\label{thm:Harm-nonzero-iff-multiple-monopoles}
  $\Harm \neq 0$ iff there are at least two linearly independent monopoles at one (equivalently, every) vertex $x$. 
\end{cor}

\begin{defn}\label{def:boundary-term-is-nonvanishing}
  The phrase ``\emph{the boundary term is nonvanishing}'' indicates that \eqref{eqn:E(u,v)=<u_0,Lapv>+sum(normals)} holds with nonzero boundary sum when applied to $\la u,v\ra_\energy$, for every representative of $u$ except one; namely, the one specified by $u(x) = \la u,\monov\ra_\energy$.
\end{defn}

Recall from Remark~\ref{rem:transient-iff-monopoles} that the network is transient iff there are monopoles in \HE. From the Discrete Gauss-Green theorem, we obtain three more criteria for transience of the random walk.

\begin{theorem}\label{thm:transience}
  The random walk on the network $(G,\cond)$ with transition probabilities $p(x,y) = \frac{\cond_{xy}}{\cond(x)}$ is transient if and only if any of the following equivalent conditions are satisfied:
  \begin{enumerate}[(i)]
    \item the boundary term is nonvanishing,
    \item $f_k := (\ge_k + \Lap)^{-1} \gd_x$ is weak-$\ad$ convergent for some sequence $\ge_k \to 0$, or
    \item $(\ran \LapM^\ad)^{c\ell} = \Fin$.
  \end{enumerate}
\end{theorem}
Note that on any network, $\opclosure{(\ran \LapM)} \ci \Fin$ and hence $\Harm \ci \ker \LapM^\ad$; cf.~\cite[Lemma~4.8]{DGG}.   

\begin{remark}\label{rem:grounded-energy}
  An alternative approach to studying the space of finite-energy functions comes by considering the \emph{grounded inner product}
  \linenopax
  \begin{align*}
    \la u, v \ra_\gdd := \cj{u(o)}v(o) + \la u, v \ra_\energy,
  \end{align*}
  which makes $\dom \energy$ into a Hilbert space \Gdd which we call the \emph{grounded energy space}. This approach is discussed in \cite{Lyons:ProbOnTrees}, \cite{Soardi94} and \cite{Kayano88,Kayano84,MuYaYo,Yamasaki79}.
  
  Let \Gddo be the closure of $\spn\{\gd_x\}$ in \Gdd. 
  If $P_{\Gddo}$ is the projection to \Gddo, and it is applied to the constant function \one then $P_{\Gddo}\one = \one$ if and only if the network is recurrent. In fact, when the network is transient, then (modulo additive constants) both $P_{\Gddo}\one$ and $P_{\Gddo}^\perp\one$ are scalar multiples of monopoles at $o$. The space $\Gddo^\perp$ is spanned by monopoles and harmonic functions.
  See \cite[\S4.1]{DGG} for more details.
\end{remark}

%% file: resistance-metrics.tex

\section{Effective resistance}
\label{sec:Effective-resistance}

There is a natural notion of distance on finite networks, which is defined in terms of resistance. Consider each edge of the network to be an electrical resistor of resistance $\cond_{xy}^{-1}$. The \emph{effective resistance metric} $R(x,y)$ is the voltage drop between the vertices $x$ and $y$ if a current of one amp is inserted into the network at $x$ and withdrawn at $y$. It is a bit surprising that this actually gives a metric, and there are several other equivalent formulations, most of which are well-known. The essential reference for effective resistance is \cite{Kig03}, but the reader may also find the excellent treatments in \cite{Soardi94} and \cite{Lyons:ProbOnTrees} to be helpful.

\begin{theorem}\label{thm:effective-resistance-metric}
  The resistance $R(x,y)$ has the following equivalent formulations:
  \linenopax
  \begin{align}
    R(x,y)
    &= \{v(x)-v(y) \suth \Lap v = \gd_x-\gd_y\}
    \label{eqn:def:R(x,y)-Lap} \\
    &= \{\energy(v) \suth \Lap v = \gd_x-\gd_y\} \label{eqn:def:R(x,y)-energy} \\
    &= 1/\min \{\energy(v) \suth v(x)=1, v(y)=0, v \in \dom\energy\} \label{eqn:def:R(x,y)-R} \\
    &= \min \{\gk \geq 0 \suth |v(x)-v(y)|^2 \leq \gk \energy(v), v \in \dom\energy\} \label{eqn:def:R(x,y)-S} \\
    &= \sup \{|v(x)-v(y)|^2 \suth \energy(v) \leq 1, v \in \dom\energy\} \label{eqn:def:R(x,y)-sup}.
  \end{align}
\end{theorem}

\begin{figure}
    \centering
    \includegraphics{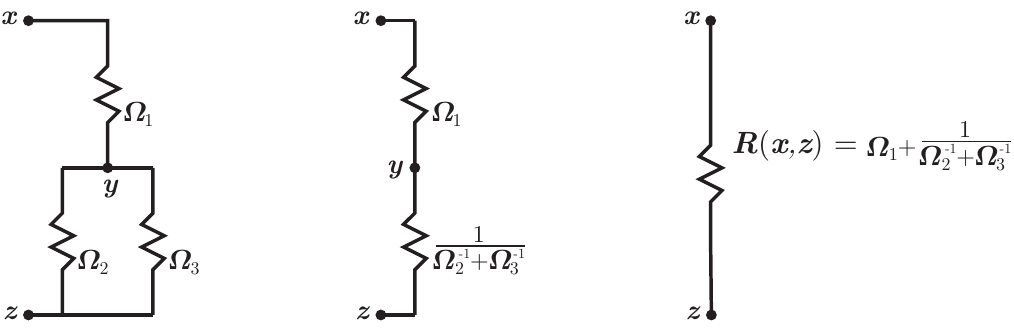}
    \caption{\captionsize Effective resistance as network reduction to a trivial network. This basic example uses parallel reduction followed by series reduction; see Remark~\ref{rem:Resistance-distance-via-network-reduction}.}
    \label{fig:network-reduction}
  \end{figure}

\begin{remark}[Resistance distance via network reduction]
  \label{rem:Resistance-distance-via-network-reduction}
  Let $G$ be a finite planar network and pick any $x,y \in \verts$. Then $G$ may be reduced to a trivial network consisting only of these two vertices and a single edge between them via the use of three basic transformations: (i) series reduction, (ii) parallel reduction, and (iii) the $\nabla$-\textsf{Y} transform \cite{Epifanov66,Truemper89}. The effective resistance between $x$ and $y$ may be interpreted as the resistance of the resulting single edge; see Figure~\ref{fig:network-reduction}.  See also \cite{Kig01} or \cite{Str06} for the $\nabla$-\textsf{Y} transform.
\end{remark}

\subsection{Resistance metric on infinite networks}
\label{sec:Resistance-metric-on-infinite-networks}

There are challenges in extending the notion of effective resistance to infinite networks. The existence of nonconstant harmonic functions $h \in \dom \energy$ implies the nonuniqueness of solutions to $\Lap u = f$ in \HE, and hence \eqref{eqn:def:R(x,y)-Lap} and \eqref{eqn:def:R(x,y)-energy} are no longer well-defined. 
This issue is studied in detail in \cite{ERM}, and in \cite{Kig03} (by very different methods). There are also accounts in \cite{Lyons:ProbOnTrees} and the literature on ``uniqueness of currents'' in infinite networks, e.g. \cite{SoardiWoess91,Thomassen90}.

Two natural choices for extension lead to the free resistance $R^F$ and the wired resistance $R^W$. In general, one has $R^F(x,y) \geq R^W(x,y)$ with equality iff $\Harm = 0$. Both of these correspond to the selection of certain solutions to $\Lap u = \gd_x - \gd_y$ (in fact, these can be interpreted as Neumann and Dirichlet boundary conditions, respectively; see \cite[Rem.~2.23]{ERM}). Also, both are given in terms of limits computed with respect to certain networks associated to an exhaustion, in the sense of Definition~\ref{def:exhaustion-of-G}. The notation $\{G_k\}_{k=1}^\iy$ always denotes an \emph{exhaustion} of the infinite network $(\Graph,\cond)$, as in Definition~\ref{def:exhaustion-of-G}. Since $x$ and $y$ are contained in all but finitely many $G_k$, we may always assume that $x,y \in G_k$, $\forall k$.

\begin{defn}\label{def:relative-resistance}
  If $H$ is a finite subnetwork of \Graph which contains $x$ and $y$, define $R_H(x,y)$ to be the \emph{resistance distance from $x$ to $y$ as computed within $H$}. In other words, compute $R_H(x,y)$ by any of the equivalent formulas of Theorem~\ref{thm:effective-resistance-metric}, but extremizing over only those functions whose support is contained in $H$.
\end{defn}

\begin{defn}\label{def:full-subnetwork}
  Let $\verts[H] \ci \verts$. Then the \emph{full subnetwork} on \verts[H] has all the edges of \Graph for which both endpoints lie in \verts[H], with the same conductances. That is, $\cond^H = \cond^G|_{\verts[H] \times \verts[H]}$.
\end{defn}

\subsubsection{Free resistance}

\begin{defn}\label{def:free-resistance}
  For any subset $\verts[H] \ci \verts$, the \emph{free subnetwork} $H^F$ is just the full subnetwork with vertices $\verts[H]$. That is, all edges of \Graph with endpoints in \verts[H] are edges of $H^F$, with the same conductances. Thus, we will denote $H^F$ by $H$ to reduce notation. 
  Let $R_{H}(x,y)$ denote the effective resistance between $x$ and $y$ as computed in $H$, as in Definition~\ref{def:relative-resistance}. 
  The \emph{free resistance} between $x$ and $y$ is defined to be
  \linenopax
  \begin{align}\label{eqn:def:free-resistance}
    R^F(x,y) := \lim_{k \to \iy} R_{G_k}(x,y),
  \end{align}
  where $\{G_k\}$ is any exhaustion of \Graph.
\end{defn}

The name ``free'' comes from the fact that this formulation is free of any boundary conditions or considerations of the complements of the $G_k$; see \cite[\S9]{Lyons:ProbOnTrees}. Theorem~\ref{thm:free-resistance} is the free extension of Theorem~\ref{thm:effective-resistance-metric} to infinite networks.

\begin{theorem}[{\cite[Thm.~2.14]{ERM}}]\label{thm:free-resistance}
  For an infinite network \Graph, the free resistance $R^F(x,y)$ has the following equivalent formulations:
  \linenopax
  \begin{align}
    R^F(x,y)
    &= v(x) - v(y), \q v=v_x-v_y  \label{eqn:def:R^F(x,y)-Lap} \\
    &= \energy(v), \q v=v_x-v_y \label{eqn:def:R^F(x,y)-energy} \\
    &= \min \{\diss(\curr) \suth \curr \in \Flo(x,y) \text{ and } \curr = \textstyle \sum \gx_\cpath \charfn{\cpath}\} \label{eqn:def:R^F(x,y)-diss} \\
    &= \left(\min\{\energy(u) \suth u \in \HE, |u(x)-u(y)|=1\}\right)^{-1} \label{eqn:def:R^F(x,y)-R} \\
    &= \inf\{\gk \geq 0 \suth |v(x)-v(y)|^2 \leq \gk \energy(v), \forall v \in \HE\} \label{eqn:def:R^F(x,y)-S} \\
    &= \sup\{|v(x)-v(y)|^2 \suth v \in \HE, \|v\|_\energy \leq 1\} \label{eqn:def:R^F(x,y)-sup} 
  \end{align}
\end{theorem}

Fix $x, y \in \Graph$ and define the operator $L_{xy}$ on \HE by $L_{xy}v := v(x)-v(y)$. 
Then \eqref{eqn:def:R^F(x,y)-S}--\eqref{eqn:def:R^F(x,y)-sup} are equivalent to $R^F(x,y) = \|L_{xy}\|$.

\subsubsection{Wired resistance}

\begin{defn}\label{def:wired-resistance}
  Given a finite full subnetwork $H$ of \Graph, define the wired subnetwork $H^W$ by identifying all vertices in $\verts \less \verts[H]$ to a single, new vertex labeled \iy. Thus, the vertex set of $H^W$ is $\verts[H] \cup \{\iy_H\}$, and the edge set of $H^W$ includes all the edges of $H$, with the same conductances. However, if $x \in \verts[H]$ has a neighbour $y \in \verts \less \verts[H]$, then $H^W$ also includes an edge from $x$ to \iy with conductance
  \linenopax
  \begin{align}\label{eqn:cond-to-iy}
    \cond_{x\iy_{\negsp[2]\scalebox{0.40}{$H$}}} := \sum_{y \nbr x, \, y \in H^\complm} \negsp[13]\cond_{xy}.
  \end{align}
  The identification of vertices in $G_k^\complm$ may result in parallel edges; then \eqref{eqn:cond-to-iy} corresponds to replacing these parallel edges by a single edge according to the usual formula for resistors in parallel.

  Let $R_{H^W}(x,y)$ denote the effective resistance between $x$ and $y$ as computed in $H^W$, as in Definition~\ref{def:relative-resistance}.
  The \emph{wired resistance} is then defined to be
  \linenopax
  \begin{align}\label{eqn:def:wired-resistance}
    R^W(x,y) := \lim_{k \to \iy} R_{G_k^W}(x,y),
  \end{align}
  where $\{G_k\}$ is any exhaustion of \Graph.
\end{defn}

  \begin{figure}
    \centering
    \includegraphics{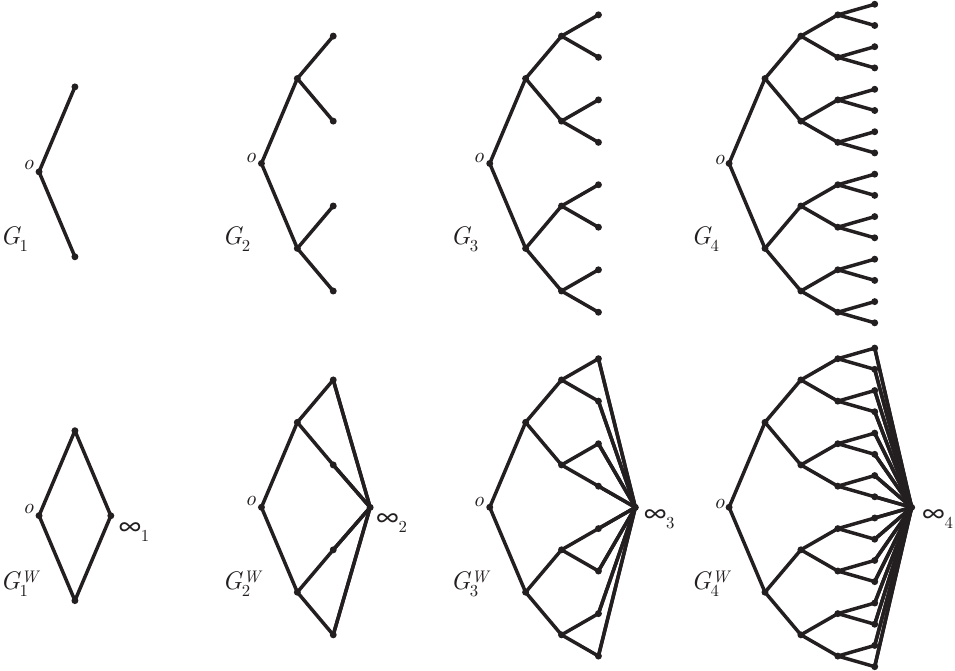}
    \caption{\captionsize Comparison of free and wired exhaustions for the example of the binary tree; see Definition~\ref{def:free-resistance} and Definition~\ref{def:wired-resistance}. Here, the vertices of $G_k$ are all those which lie within $k$ edges (``steps'') of the origin. If the edges of \Graph all have conductance 1, then so do all the edges of each $G_k^F$ and $G_k^W$, except for the edges incident upon $\iy_k = \iy_{G_k}$, which have conductance 2.}
    \label{tree-exhaustions}
  \end{figure}

The wired subnetwork is equivalently obtained by ``shorting together'' all vertices of $H^\complm$, and hence it follows from Rayleigh's monotonicity principle that $R^W(x,y) \leq R^F(x,y)$; cf. \cite[\S1.4]{DoSn84} or \cite[\S2.4]{Lyons:ProbOnTrees}.

\begin{theorem}[{\cite[Thm.~2.20]{ERM}}]\label{thm:wired-resistance}
  The wired resistance may be computed by any of the following equivalent formulations:
  \linenopax
  \begin{align}
    R^W(x,y)
    &= f(x)-f(y), \q f = f_x - f_y  \label{eqn:def:R^W(x,y)-Lap} \\
    &= \energy(f), \q f = f_x - f_y \label{eqn:def:R^W(x,y)-energy} \\
    &= \smash{\left(\min\{\energy(v) \suth |v(x) - v(y)|=1, v \in \Fin\}\right)^{-1}} \label{eqn:def:R^W(x,y)-R} \\
    &= \inf\{\gk \geq 0 \suth |v(x)-v(y)|^2 \leq \gk \energy(v), \forall v \in \Fin\} \label{eqn:def:R^W(x,y)-S} \\
    &= \sup\{|v(x)-v(y)|^2 \suth v \in \Fin, \|v\|_\energy \leq 1\} \label{eqn:def:R^W(x,y)-sup} 
  \end{align}
\end{theorem}
 
Note that \eqref{eqn:def:R^W(x,y)-Lap} and \eqref{eqn:def:R^W(x,y)-energy} are  equivalent to 
\begin{align*}
  R^W(x,y) 
  &= \min \{v(x)-v(y) \suth \Lap v = \gd_x-\gd_y, v \in \dom\energy\} \tag{\ref{eqn:def:R^W(x,y)-Lap}'}\\
  &= \min \{\energy(v) \suth \Lap v = \gd_x-\gd_y, v \in \dom\energy\} \tag{\ref{eqn:def:R^W(x,y)-energy}'}.
\end{align*}

\subsection{von Neumann construction of the energy space \HE}
\label{sec:vonNeumann's-embedding-thm}

Let $R=R^F$ or $R=R^W$. The discussion of the effective resistance is important in this paper in two respects. 
\begin{enumerate}[(i)]
  \item Theorem~\ref{thm:R^F-embed-ERN-in-Hilbert} shows that \HE is the natural Hilbert space for studying the metric space $(\Graph,R)$.
  \item The function $R(x,y)$ allows us to construct a probability measure in Theorem~\ref{thm:HE-isom-to-L2(S',P)}.
\end{enumerate}
Both of these results stem from the fact that (free or wired) effective resistance is a negative semidefinite function on $\verts \times \verts$, as is shown in \cite[Thm.~5.4]{ERM}.

\begin{defn}\label{def:negative-semidefinite}
  A function $M:X \times X \to \bR$ is \emph{negative semidefinite} iff for any $f:X \to \bR$ satisfying $\sum_{x \in X} f(x) = 0$, one has
  \linenopax
  \begin{equation}\label{eqn:negative-semidefinite}
    \sum_{x,y \in F} f(x) M(x,y) f(y) \leq 0,
  \end{equation}
  where $F$ is any finite subset of $X$. 
\end{defn}

One can think of $M$ as a matrix and \eqref{eqn:negative-semidefinite} as matrix multiplication. von Neumann and Schoenberg \cite{vN32a,Ber96,Schoe38a,Schoe38b,Ber84} showed that \eqref{eqn:negative-semidefinite} is precisely the condition that allows one to embed a metric space into a Hilbert space. This theorem also has a form of uniqueness which may be thought of as a universal property.

\begin{theorem}[von Neumann]\label{thm:vonNeumann's-embedding-thm}
  Suppose $(X,d)$ is a metric space. There exists a Hilbert space $\sH$ and an embedding $w:(X,d) \to \sH$ sending $x \mapsto w_x$ and satisfying
  \linenopax
  \begin{equation}\label{eqn:vonNeumann's-embedding-thm}
    d(x,y) = \|w_x - w_y\|_\sH
  \end{equation}
  if and only if $d^2$ is negative semidefinite. 
  
  Furthermore, if there is another Hilbert space \sK and an embedding $k:\sH \to \sK$, with $\|k_x-k_y\|_\sK = d(x,y)$ and $\{k_x\}_{x \in X}$ dense in \sK, then there exists a unique unitary isomorphism $U:\sH \to \sK$.
\end{theorem}

\begin{theorem}[{\cite[Thm.~5.4]{ERM}}]\label{thm:R^F-embed-ERN-in-Hilbert}
  $(\Graph,R^F)$ may be isometrically embedded in a Hilbert space, and this Hilbert space is unitarily equivalent to \HE. Under this embedding, $x$ is mapped to the energy kernel element $v_x$.
  
  Moreover, $(\Graph,R^W)$ may be isometrically embedded in a Hilbert space, and this Hilbert space is unitarily equivalent to \Fin. Under this embedding, $x$ is mapped to $f_x = \Pfin v_x$.
\end{theorem}

By this theorem, we see that \HE is the natural choice of Hilbert space for studying the metric spaces $(\Graph, R^F)$ and $(\Graph, R^W)$.

%% file: gelfand.tex

\section{A boundary integral representation for the harmonic functions}
\label{sec:boundary-integral-representation-for-harmonics}

We are motivated by the following result, which follows readily from Theorem~\ref{thm:E(u,v)=<u,Lapv>+sum(normals)} and may be found in \cite[Cor.~3.14]{DGG}.

\begin{theorem}[Boundary representation of harmonic functions]
  \label{thm:Boundary-representation-of-harmonic-functions}
  For $u \in \spn\{h_x\}$, 
  \linenopax
  \begin{align}\label{eqn:Boundary-representation-of-harmonic-functions}
    u(x) = \sum_{\bd \Graph} u \dn{h_x} + u(o). 
  \end{align}
  \begin{proof}
    Note that 
    $u(x) - u(o) 
    = \la v_x, u\ra_\energy 
    = \cj{\la u, v_x\ra_\energy} 
    = \sum_{\bd \Graph} u \dn{h_x}$ by \eqref{eqn:def:vx}.
  \end{proof}
\end{theorem}

Formula \eqref{eqn:Boundary-representation-of-harmonic-functions} begs comparison with the Poisson integral formula.
Recall the classical result of Poisson that gives a kernel $k:\gW \times \del \gW \to \bR$ from which a bounded harmonic function can be given via
\linenopax
  \begin{equation}\label{eqn:Poisson-bdy-repn}
    u(x) = \int_{\del \gW} u(y) k(x,dy),
    \qq y \in \del \gW.
  \end{equation}

One would like to obtain a (probability) measure space to serve as the boundary of \Graph. We use some techniques from the theory of stochastic integration for which it was shown in \cite{Nelson64} that a Hilbert space does not suffice; see \cite[\S3.1]{Hida80}. The workaround {is} to build a \emph{Gel'fand triple} $S \ci \sH \ci S'$ (a more precise definition appears just below), and construct a suitable probability measure on $S'$. In \S\ref{sec:Gel'fand-triples-and-duality}, we briefly describe the general theory of Gel'fand triples as they apply in the current context. In \S\ref{sec:Gel'fand-triple-for-HE}, we use \Lap to construct a Gel'fand triple for \HE. Then in \S\ref{sec:Wiener-embedding}, we apply the general theory to the Gel'fand triple $\Schw \ci \HE \ci \Schw'$ and obtain a Gaussian probability measure \prob on $\Schw'$, and an isometric embedding $\HE \hookrightarrow L^2(\Schw',\prob)$. This allows us to study the boundary $\bd G$ as a subset of $\Schw'$. 
For the general theory of analysis in Hilbert space, see \cite{Gro67,Gro70}.

\subsection{Gel'fand triples and duality}
\label{sec:Gel'fand-triples-and-duality}

In a little more detail, a \emph{Gel'fand triple} (also called a \emph{rigged Hilbert space}) is
\linenopax
\begin{equation}\label{eqn:Gel'fand-triple-intro}
  S \ci \sH \ci S',
\end{equation}
where $S$ is dense in \sH and $S'$ is the dual of $S$. 
While $S$ is a dense subspace of \sH with respect to the Hilbert norm, it also comes equipped with a strictly finer ``test function'' topology, and it is required that the inclusion mapping of $S$ into \sH is continuous with respect to these topologies. Therefore, when the dual $S'$ is taken with respect to this finer topology, one obtains a \emph{strict} containment $\HE \subsetneq S'$. It turns out that $S'$ is large enough to support a (Gaussian!) probability measure. 

We will give a ``test function topology'' as a Fr\'{e}chet topology defined via a specific sequence of seminorms. It was Gel'fand's idea to formalize this construction abstractly using a system of nuclearity axioms \cite{GMS58, Minlos58, Minlos59}. This presentation is adapted from quantum mechanics. 

\begin{remark}[Tempered distributions and the Laplacian]
  \label{rem:tempered-distributions}
  There is a concrete situation when the Gel'fand triple construction is especially natural: $\sH = L^2(\bR,dx)$ and $S$ is the \emph{Schwartz space} of functions of rapid decay. That is, each $f \in S$ is a $C^\iy$ smooth function which decays (along with all its derivatives) faster than any polynomial as $x \to \pm \iy$. In this case, $S'$ is the space of \emph{tempered distributions} and the seminorms defining the Fr\'{e}chet topology on $S$ are 
\linenopax
\begin{align*}
  p_m(f) := \sup \{|x^k f^{(n)}(x)| \suth x \in \bR, 0 \leq k,n \leq m\},
  \qq m=0,1,2,\dots,
\end{align*}
where $f^{(n)}$ is the \nth derivative of $f$. Then $S'$ is the dual of $S$ with respect to this Fr\'{e}chet topology. One can equivalently express $S$ as
\linenopax
\begin{align}\label{eqn:Schwartz-space-as-powers-of-Hamiltonian}
  S := \{f \in L^2(\bR) \suth (\tilde{P}^2 + \tilde{Q}^2)^n f \in L^2(\bR), \forall n\},
\end{align}
where $\tilde{P}:f(x) \mapsto \frac1\ii\frac{d}{dx}$ and $\tilde{Q}:f(x) \mapsto x f(x)$ are Heisenberg's operators. The operator $\tilde{P}^2 + \tilde{Q}^2$ is often called the quantum mechanical Hamiltonian, but some others (e.g., Hida, Gross) would call it a Laplacian, and this perspective tightens the analogy with the present context. In this sense, \eqref{eqn:Schwartz-space-as-powers-of-Hamiltonian} could be rewritten $S := \dom \Lap^\iy$; compare to \eqref{eqn:Schwartz} just below.
\end{remark}

  The duality between $S$ and $S'$ allows for the extension of the inner product on \sH to a pairing of $S$ and $S'$:
\linenopax
\begin{align}\label{eqn:extended-pairing}
  \la \cdot,\cdot\ra_\sH:\sH \times \sH \to \bC
  \qq\text{to}\qq
  \la \cdot,\cdot\ra_{\tilde \sH}: S \times S' \to \bR.
\end{align}
In other words, one obtains a Fourier-type duality restricted to $S$. 

The proof of Theorem~\ref{thm:HE-isom-to-L2(S',P)} will require Minlos' generalization of Bochner's theorem from \cite{Minlos63, Sch73}. This important result states that a cylindrical measure on the dual of a nuclear space is a Radon measure iff its Fourier transform is continuous. In this context, however, the notion of Fourier transform is infinite-dimensional; cf.~\cite{Lee96}. 

\begin{theorem}[Minlos]
  \label{thm:Minlos'-theorem}
  Given a Gel'fand triple $S \ci \sH \ci S'$, 
  there is a bijective correspondence between the positive definite functions $f$ on $S$ and the Radon probability measures on $S'$, determined uniquely by the identity
  \linenopax
  \begin{align}\label{eqn:Minlos-identity}
    f(s) = \int_{S'} e^{\ii \la s, \gx\ra_{\tilde \sH}} \,d\prob_{\negsp[6]f}(\gx),
    \qq\forall s \in S,
  \end{align}
  where 
  $\la \cdot , \cdot \ra_{\tilde \sH}$ is the extended pairing on $S \times S'$ as in \eqref{eqn:extended-pairing}. 
\end{theorem}

Formula \eqref{eqn:Minlos-identity} may be interpreted as defining the Fourier transform of \prob. We apply Minlos' theorem in the standard manner for white noise constructions, and obtain the following corollary.

\begin{cor}[White noise]
  \label{thm:White-noise-formula}
  Given a Gel'fand triple $S \ci \sH \ci S'$, there is a probability measure \prob on $S'$ satisfying
  \begin{align}\label{eqn:White-noise-formula}
    e^{-\frac12 \la s,s\ra_\sH}
    =\int_{S'} e^{\ii \la s, \gx\ra_{\tilde \sH}} \,d\prob(\gx).
  \end{align}
\end{cor}

In the proof of Theorem~\ref{thm:HE-isom-to-L2(S',P)}, we show that \prob in \eqref{eqn:White-noise-formula} is actually a Gaussian measure on $\Schw'$.
The function on the left-hand side of \eqref{eqn:White-noise-formula} 
plays a special role in stochastic integration, and its use in quantization. To see that it is a positive definite function on $S$, we appeal to a famous result of Schoenberg which may be found in \cite{Ber84,ScWh49}.

\begin{theorem}[Schoenberg]
  \label{thm:Schoenberg's-Thm}
  Let $X$ be a set and let $Q: X \times X \to \bR$ be a function. Then the following are equivalent.
  \begin{enumerate}
    \item $Q$ is negative semidefinite.
    \item $\forall t \in \bR^+$, the function $p_t(x,y) := e^{-t Q(x,y)}$ is positive definite on $X \times X$.
    \item There exists a Hilbert space $\sH$ and a function $f:X \to \sH$ such that \[Q(x,y) = \|f(x)-f(y)\|_\sH^2.\]
  \end{enumerate}
\end{theorem}

In the proof of Theorem~\ref{thm:HE-isom-to-L2(S',P)}, we apply Schoenberg's Theorem with $t=\frac12$ to the resistance metric in the form 
\linenopax
\begin{align}\label{eqn:def:R^F(x,y)-energy}
  R^F(x,y) = \|v_x-v_y\|_\energy^2,
\end{align}
which appears in \cite[Thm.~2.13]{ERM}. Recall from \S\ref{sec:vonNeumann's-embedding-thm} that \eqref{eqn:def:R^F(x,y)-energy} is negative semidefinite.

\subsection{A Gel'fand triple for \HE}
\label{sec:Gel'fand-triple-for-HE}

To apply Minlos' Theorem, we first need to construct a Gel'fand triple for \HE; we begin by identifying a certain subspace of $\MP = \dom \LapM$ (as given in Definition~\ref{def:LapM}) to act as the space of test functions, which we denote \Schw.

\begin{defn}\label{def:extn-of-Lap}
  Let \LapS be a self-adjoint extension of \LapM; since \LapM is Hermitian and commutes with conjugation (since \cond is \bR-valued), a theorem of von Neumann's states that such an extension exists.
\glossary{name={\LapS},description={a self-adjoint extension of \LapM},sort=L,format=textbf}
  
  Let $\LapS^p u := (\LapS\LapS\dots\LapS)u$ be the $p$-fold product of \LapS applied to $u \in \HE$. Define $\dom(\LapS^p)$ inductively by
  \linenopax
  \begin{equation}\label{eqn:domE(LapEp)}
    \dom(\LapS^p) := \{u \suth \LapS^{p-1} u \in \dom(\LapS)\}.
  \end{equation}
\end{defn}

\begin{defn}[Test functions]
  \label{def:Schw-on-G}
  The \emph{(Schwartz) space of potentials of rapid decay} is 
  \linenopax
  \begin{equation}\label{eqn:Schwartz}
    \Schw := \dom(\LapS^\iy), 
  \end{equation}
\glossary{name={$\Schw$},description={``Schwartz space'' of test functions (of rapid decay)},sort=S,format=textbf}
  where $\dom(\LapS^\iy) := \bigcap_{p=1}^\iy \dom(\LapS^p)$ consists of all $u \in \HE$ for which $\LapS^p u \in \HE$ for any $p$. 
\end{defn}

\begin{defn}[Distributions]
  \label{def:Frechet-topology}\label{def:Schwartz-distributions}
  For each $p \in \bN$, there is a seminorm on \Schw defined by
  \linenopax
  \begin{equation}\label{eqn:p-norm-on-Schw}
    \|u\|_p := \| \LapS^p u\|_\energy.
  \end{equation}
Since $(\dom \LapS^p, \|\cdot\|_p)$ is a Hilbert space for each $p \in \bN$, the system of seminorms $\sP = \{\|\cdot\|_p\}_{p \in \bN}$ defines a Fr\'{e}chet topology on \Schw. 
   The space $\Schw'$ of \emph{Schwartz distributions} or \emph{tempered distributions} is the (dual) space of \sP-continuous linear functionals on \Schw.
\end{defn}

\begin{remark}\label{rem:Schw-contains-vx}
  If $\deg(x)$ is finite for each $x \in \verts$, or if $\uBd < \iy$, then one has $v_x \in \Schw$. In the first case, this can be proved from the identity $\gd_x = \cond(x) v_x - \sum_{y \nbr x} \cond_{xy} v_y$ which is given in \cite[Lem.~2.22]{DGG}. In the second case, the bound on \cond implies \LapS is bounded and hence everywhere-defined.

  \version{}{\marginpar{cut this?}}
  When \Schw contains $\{v_x\}$, it should be noted that $\spn\{v_x\}$ is dense in \Schw with respect to \energy, but \emph{not} with respect to the Fr\'{e}chet topology induced by the seminorms \eqref{eqn:p-norm-on-Schw}, nor with respect to the graph norm. One has the inclusions
  \linenopax
  \begin{equation}\label{eqn:vx-graph-inclusions}
    \left\{\left[\begin{array}{c} v_x \\ \LapM v_x \end{array}\right]\right\}
    \ci \left\{\left[\begin{array}{c} s \\ \LapS s \end{array}\right]\right\}
    \ci \left\{\left[\begin{array}{c} u \\ \LapS u \end{array}\right]\right\}
  \end{equation}
  where $s \in \Schw$ and $u \in \HE$. The second inclusion is dense but the first is not.
\end{remark}

\begin{remark}\label{rem:Schwartz-space-is-real-valued}
  Note that \Schw and $\Schw'$ consist of \bR-valued functions. This technical detail is important because we do not expect the integral $\int_{S'} e^{\ii \la u, \cdot \ra_{\tilde \Wiener}} \,d\prob$ from \eqref{eqn:Minlos-identity} to converge unless it is certain that $\la u, \cdot \ra$ is \bR-valued. This is the reason for the last conclusion of Theorem~\ref{thm:Wiener-product-as-Lap(p)-powers}.
\end{remark}
    
\begin{defn}\label{def:spectral-truncation}
  Let $\charfn{[a,b]}$ denote the usual indicator function of the interval $[a,b] \ci \bR$, and let \spectrans be the spectral transform in the spectral representation of \LapS, and let $E$ be the associated projection-valued measure. Then define $E_n$ to be the \emph{spectral truncation operator} acting on \HE by
  \linenopax
  \begin{align*}
     E_n u 
     := \spectrans^\ad \charfn{[\frac1n,n]} \spectrans u 
     = \int_{1/n}^n E(dt)u.
  \end{align*}
\end{defn}


\begin{lemma}\label{thm:Schw-dense-in-HE}
  \version{}{\marginpar{Cut this?}}
  With respect to \energy, \Schw is a dense analytic subspace of \HE.
  \begin{proof}
    This essentially follows immediately once it is clear that $E_n$ maps \HE into \Schw. For $u \in \HE$, and for any $p=1,2,\dots$,
    \linenopax
    \begin{equation}\label{eqn:E-norm-of-Lapp(spectral-truncation)}
      \|\LapS^{p} E_n u\|_\energy^2 
      = \int_{1/n}^n \gl^{2p} \|E(d\gl)u\|_\energy^2
      \leq n^{2p} \|u\|_\energy^2,
    \end{equation}
    So $E_n u \in \Schw$. It follows that $\|u-E_n u\|_\energy \to 0$ by standard spectral theory.
  \end{proof}
\end{lemma}

\begin{theorem}\label{thm:Wiener-product-as-Lap(p)-powers}
  $\Schw \ci \HE \ci \Schw'$ is a Gel'fand triple, and the energy form $\la \cdot,\cdot\ra_\energy$ extends to a pairing on $\Schw \times \Schw'$ defined by
 \linenopax
  \begin{equation}\label{eqn:energy-extends-by-ps}
    \la u, \gx\ra_\Wiener := \la \LapS^p u, \LapS^{-p} \gx \ra_\energy,
  \end{equation}
  where $p$ is any integer such that $| \gx(u)| \leq K \|\Lap^p u\|_\energy$ for all $u \in \Schw$. This pairing on $\Schw \times \Schw'$ is equivalently given by
  \linenopax
  \begin{equation}\label{eqn:energy-extends-by-vn}
    \la u,\gx\ra_\Wiener = \lim_{n \to \iy} \gx(E_n u), 
  \end{equation}
  where the limit is taken in the topology of $\Schw'$. 
\end{theorem}

\begin{cor}\label{thm:Spectruncation-extends-to-S'}
  $E_n$ extends to a mapping $\tilde E_n: \Schw' \to \HE$ defined via $\la u, \tilde E_n \gx \ra_\energy := \gx(E_n u)$. Thus, we have a pointwise extension of $\la \cdot \,,\,\cdot\ra_\Wiener$ to $\HE \times \Schw'$ given by
  \begin{equation}\label{eqn:energy-extends-by-truncating-on-S'}
    \la u,\gx\ra_\Wiener = \lim_{n \to \iy} \la u, \tilde E_n \gx \ra_\energy.
  \end{equation}
\end{cor}

%% file: wiener.tex

\subsection{The Wiener embedding and the space $\Schw'$}
\label{sec:Wiener-embedding}

With Theorem~\ref{thm:Wiener-product-as-Lap(p)-powers}, we have a Gel'fand triple and we are now ready to apply the white noise construction of Cor.~\ref{thm:White-noise-formula}. 
%
Note that in Theorem~\ref{thm:HE-isom-to-L2(S',P)}, expectations are taken with respect to the variable $\gx \in \Schw'$, that is, $\Ex(f) := \int_{\Schw'} f(\gx) \,d\prob(\gx)$. 
\glossary{name={$\Ex(f)$},description={expectation of $f$ with dummy variable \gx; $\int f(\gx)\,d\gm(x)$},sort=E,format=textbf}

\begin{theorem}[Wiener embedding]
  \label{thm:HE-isom-to-L2(S',P)}
  The Wiener transform $\sW:\HE \to L^2(\Schw',\prob)$ is given by
  \linenopax
  \begin{equation}\label{eqn:Gaussian-transform}
    \sW : v \mapsto \tilde v,
    \q \tilde v(\gx) := \la v, \gx\ra_\Wiener,
  \end{equation}
  and is an isometry. The extended reproducing kernel $\{\tilde v_x\}_{x \in \verts}$ is a system of Gaussian random variables which gives the resistance distance by
  \linenopax
  \begin{equation}\label{eqn:R(x,y)-as-expectation}
    R^F(x,y) = \Ex((\tilde v_x - \tilde v_y)^2).
  \end{equation}
  Moreover, for any $u,v \in \HE$, the energy inner product extends directly as
  \linenopax
  \begin{equation}\label{eqn:Expectation-formula-for-energy-prod}
    \la u, v \ra_\energy
    = \Ex\left( \cj{\tilde{u}} \tilde{v} \right)
    = \int_{\Schw'} \cj{\tilde{u}} \tilde{v} \,d\prob.
  \end{equation}
  \begin{proof}
    Since $R^F(x,y)$ is negative semidefinite (see \cite[Thm.~5.4]{ERM}), we may apply Schoenberg's theorem and deduce that $\exp(-\tfrac12\|u-v\|_\energy^2)$ is a positive definite function on $\HE \times \HE$. Consequently, an application of the Minlos correspondence (Theorem~\ref{thm:Minlos'-theorem}) to the Gel'fand triple established in Lemma~\ref{thm:Schw-dense-in-HE} yields a Gaussian probability measure \prob on $\Schw'$. 
    
    Moreover, \eqref{eqn:Minlos-identity} gives 
    \linenopax
    \begin{align}\label{eqn:Minlos-eqns}
      \Ex(e^{\ii\la u, \gx \ra_\Wiener}) = e^{-\frac12\|u\|_\energy^2},
    \end{align}
    provided that \gx is \bR-valued (so that the integral converges). Therefore, we give the proof for the \bR-valued subspace of \Schw (and of $\Schw'$), and then complexify in the last step via the standard decomposition into real and complex parts: $u = u_1 + \ii u_2$ where $u_i$ is a \bR-valued elements of \HE, etc.

    From \eqref{eqn:Minlos-eqns}, one computes
    \linenopax
    \begin{align}\label{eqn:Minlos-expectation-integral}
      \int_{\Schw'} \left(1 + \ii\la u, \gx\ra_\Wiener - \frac12\la u, \gx\ra_\Wiener^2 + \cdots \right)\,d\prob(\gx) 
      = 1 - \frac12 \la u, u \ra_\energy + \cdots.
    \end{align}
    Now it follows that $\Ex(\tilde{u}^2) = \Ex(\la u, \gx\ra_\Wiener^2) = \|u\|_\energy^2$ for every $u \in \Schw$, by comparing the terms of \eqref{eqn:Minlos-expectation-integral} which are quadratic in $u$. Therefore, $\sW:\HE \to \Schw'$ is an isometry, and \eqref{eqn:Minlos-expectation-integral} gives
    \linenopax
    \begin{align}\label{eqn:Exp-tilde-vx=E(vx)}
      \Ex(|\tilde v_x - \tilde v_y|^2)
      = \Ex(\la v_x - v_y, \gx \ra^2)
      &= \|v_x - v_y\|_\energy^2,
    \end{align}
    whence \eqref{eqn:R(x,y)-as-expectation} follows from \eqref{eqn:def:R^F(x,y)-energy}. Note that by comparing the linear terms, \eqref{eqn:Minlos-expectation-integral} implies $\Ex(1) = 1$, so that \prob is a probability measure, and $\Ex(\la u,\gx\ra) = 0$ and $\Ex(\la u,\gx\ra^2) = \|u\|_\Wiener^2$, so that \prob is actually Gaussian.
    
    
    Finally, use polarization to compute
    \linenopax
    \begin{align*}
      \la u, v \ra_\energy
      &= \frac14 \left(\|u+v\|_\energy^2 - \|u-v\|_\energy^2\right) \\
      &= \frac14 \left(\Ex\left(\left|\tilde{u}+ \tilde{v}\right|^2\right) 
        - \Ex\left(\left|\tilde{u}-\tilde{v}\right|^2\right)\right) 
        &&\text{by \eqref{eqn:Exp-tilde-vx=E(vx)}} \\
      &= \frac14 \int_{\Schw'} \left|\tilde{u}+\tilde{v}\right|^2(\gx)
        - \left|\tilde{u}-\tilde{v}\right|^2(\gx) \,d\prob(\gx) \\
      &= \int_{\Schw'} \cj{\tilde{u}}(\gx) \tilde{v}(\gx) \,d\prob(\gx).
    \end{align*}
    This establishes \eqref{eqn:Expectation-formula-for-energy-prod} and, upon complexification, completes the proof.
  \end{proof}
\end{theorem}
 
\begin{remark}\label{rem:wired-Wiener}
    Observe that Theorem~\ref{thm:HE-isom-to-L2(S',P)} was carried out for the free resistance, but all the arguments go through equally well for the wired resistance; note that $R^W$ is similarly negative semidefinite by Theorem~\ref{thm:Schoenberg's-Thm} and \cite[Cor.~5.5]{ERM}. Thus, there is a corresponding Wiener transform $\sW:\Fin \to L^2(\Schw',\prob)$ defined by
  \linenopax
  \begin{equation}\label{eqn:wired-Gaussian-transform}
    \sW : v \mapsto \tilde f,
    \qq f = \Pfin v \;\text{ and }\; \tilde f(\gx) = \la f, \gx\ra_\Wiener.
  \end{equation}
  Again, $\{\tilde f_x\}_{x \in \verts}$ is a system of Gaussian random variables which gives the wired resistance distance by $R^W(x,y) = \Ex((\tilde f_x - \tilde f_y)^2)$. 
\end{remark}
  
\begin{remark}
  \label{rem:abuse-of-extension-notation}
  For $u \in \Harm$ and $\gx \in \Schw'$, let us abuse notation and write $u$ for $\tilde{u}$ so as to avoid unnecessary tildes. That is, $u(\gx) := \tilde{u}(\gx) = \la u, \gx\ra_\Wiener$. 
\end{remark}

\begin{remark}\label{rem:Wiener-improves-Minlos}
  \version{}{\marginpar{Cut this?}}
  The polynomials are dense in $L^2(\Schw',\prob)$: let $\gf(t_1, t_2, \dots, t_k)$ denote an ordinary polynomial in $k$ variables. Then 
  \linenopax
  \begin{align}\label{eqn:polynomials-in-S'}
  \gf(\gx) := \gf\left({u_1}(\gx), {u_2}(\gx), \dots \vstr[2.2] \, {u_n}(\gx)\right)
  \end{align}
  is a polynomial on $\Schw'$ and 
  \linenopax
  \begin{align}\label{eqn:Poly(n)-in-S'}
    \Poly_n := \{\gf\left({u_1}(\gx),{u_2}(\gx), \dots \vstr[2.2] \,{u_k}(\gx)\right), \deg(\gf) \leq n, \suth u_j \in \HE, \gx \in \Schw'\}
  \end{align}
  is the collection of polynomials of degree at most $n$, and $\{\Poly_n\}_{n=0}^\iy$ is an increasing family whose union is all of $\Schw'$. One can see that the monomials $\la u, \gx\ra_\Wiener$ are in $L^2(\Schw',\prob)$ as follows: compare like powers of $u$ from either side of \eqref{eqn:Minlos-expectation-integral} to see that $\Ex\left(\la u, \gx\ra_\Wiener^{2n+1}\right) = 0$ and
  \linenopax
  \begin{align}\label{eqn:expectation-of-even-monomials}
    \Ex\left(\la u, \gx\ra_\Wiener^{2n}\right) 
    = \int_{\Schw'} |\la u, \gx\ra_\Wiener|^{2n} \, d\prob(\gx) 
    = \frac{(2n)!}{2^n n!} \|u\|_\energy^{2n}, 
  \end{align}
  and then apply the Schwarz inequality. 
  
  To see why the polynomials $\{\Poly_n\}_{n=0}^\iy$ should be dense in $L^2(\Schw',\prob)$ observe that the sequence $\{P_{\Poly_n}\}_{n=0}^\iy$ of orthogonal projections increases to the identity, and therefore, $\{P_{\Poly_n} \tilde u\}$ forms a martingale, for any $u \in \HE$ (i.e., for any $\tilde u \in L^2(\Schw',\prob)$).
  
  Denote the ``multiple Wiener integral of degree $n$'' by 
  \linenopax
  \begin{align*}
    H_n := 
    \left(cl \spn\{\la u, \cdot\ra_\Wiener^n \suth u \in \HE\}\right) \ominus \{\la u, \cdot\ra_\Wiener^k \suth k<n, u \in \HE\},
  \end{align*}
  for each $n \geq 1$, and $H_0 := \bC \one$ for a vector \one with $\|\one\|_2=1$. Then we have an orthogonal decomposition of the Hilbert space     
  \linenopax
  \begin{align}\label{eqn:Fock-repn-of-L2(S',P)}
    L^2(\Schw',\prob) = \bigoplus_{n=0}^\iy H_n.
  \end{align}
  See \cite[Thm.~4.1]{Hida80} for a more extensive discussion.
  A physicist would call \eqref{eqn:Fock-repn-of-L2(S',P)} the Fock space representation of $L^2(\Schw',\prob)$ with ``vacuum vector'' \one. Note that $H_n$ has a natural (symmetric) tensor product structure: $H_n \cong \HE^{\otimes n}$, the $n$-fold symmetric tensor product of \HE with itself. Observe that \one is orthogonal to \Fin and \Harm, but is not the zero element of $L^2(\Schw',\prob)$.

  Familiarity with these ideas is not necessary for the sequel, but the decomposition \eqref{eqn:Fock-repn-of-L2(S',P)} is helpful for understanding two key things:
  \begin{enumerate}[(i)]
    \item The Wiener isometry $\sW:\HE \to L^2(\Schw',\prob)$ identifies \HE with the subspace $H_1$ of $L^2(\Schw',\prob)$, in particular, $L^2(\Schw',\prob)$ is not isomorphic to \HE. In fact, it is the second quantization of \HE.
    \item The constant function \one is an element of $L^2(\Schw',\prob)$ but does not correspond to any element of \HE. In particular, \one is not equivalent to 0 in $L^2(\Schw',\prob)$ (as it was in \HE).
  \end{enumerate}
  It is somewhat ironic that we began this story by removing the constants (via the introduction of \HE), only to reintroduce them with a certain amount of effort, much later. 
\end{remark}

%% file: boundary.tex


Recall that we began with a comparison of the Poisson boundary representation for bounded harmonic functions with the boundary sum representation recalled in Theorem~\ref{thm:Boundary-representation-of-harmonic-functions}:
\linenopax
\begin{equation*}
  u(x) = \int_{\del \gW} u(y) k(x,dy) 
  \qq\leftrightarrow\qq
  u(x) = \sum_{\bd \Graph} u \dn{h_x} + u(o).
\end{equation*}
In this section, we replace the sum with an integral and complete the parallel.

\begin{cor}[Boundary integral representation for harmonic functions]
  \label{thm:Boundary-integral-repn-for-harm} \hfill \\
  For any $u \in \Harm$ and with $h_x = \Phar v_x$, 
  \linenopax
  \begin{equation}\label{eqn:integral-boundary-repn-of-h}
    u(x) = \int_{\Squoth} u(\gx) h_x(\gx) \, d\quotprob(\gx) + u(o).
  \end{equation}
  \begin{proof}
    Starting with \eqref{eqn:def:vx}, compute
    \linenopax
    \begin{align}\label{eqn:boundary-repn-for-harmonic-integral-derived}
      u(x) - u(o)
      = \la h_x, u\ra_\energy
      = \cj{\la u, h_x \ra_\energy}
      = \cj{\int_{\Schw'} \cj{u} h_x \, d\quotprob},
    \end{align}
    where the last equality comes by substituting $v=h_x$ in \eqref{eqn:Expectation-formula-for-energy-prod}. It is shown in \cite[Lem.~2.24]{DGG} that $\cj{h_x}=h_x$. 
  \end{proof}
\end{cor}

\begin{remark}[A Hilbert space interpretation of {bd}\,\Graph]
  \label{rem:boundary-of-G-from-Minlos}
  In view of Corollary~\ref{thm:Boundary-integral-repn-for-harm}, we are now able to ``catch'' the boundary between \Schw and $\Schw'$: 
  the boundary $\bd \Graph$ may be thought of as (a presumably proper subset of) $\Squoth / \HE$. In parallel to the construction of the Martin boundary, one expects that $\Squoth / \HE$ is larger than necessary, and that \prob is probably supported on a much smaller set, comparable to the minimal Martin boundary; cf.~\cite[Ch.~IV]{Woess00}. 
  Corollary~\ref{thm:Boundary-integral-repn-for-harm} suggests that $\mathbbm{k}(x,d\gx) := h_x(\gx) d\quotprob$ is the discrete analogue in \HE of the Poisson kernel $k(x,dy)$, and comparison of \eqref{eqn:boundary-repn-for-harmonic} with \eqref{eqn:integral-boundary-repn-of-h} gives a way of understanding a boundary integral as a limit of Riemann sums:
  \linenopax
  \begin{equation}\label{eqn:boundary-of-G-from-Minlos}
    \int_{\Schw'} u \, h_x \, d\quotprob
    = \lim_{k \to \iy} \sum_{\bd G_k} u(x) \dn{h_x}(x). 
  \end{equation}
  (We continue to omit the tildes as in Remark~\ref{rem:abuse-of-extension-notation}.) By a theorem of Nelson, \quotprob is fully supported on those functions which are H\"{o}lder-continuous with exponent $\ga=\frac12$, which we denote by $\Lip(\frac12) \ci \Schw'$; see \cite{Nelson64,Nel69}. Recall from \cite[Cor.~2.16]{ERM} that $\HE \ci Lip(\frac12)$.
  %
  Current research focuses on determining the precise relationship between $\bd \Graph$ and these other spaces (Martin boundary, $Lip(\frac12)$), and an explicit representation of $\bd \Graph$ in terms of paths in \Graph and/or cocycles. We expect that $\bd\Graph$ will have applications in the analysis of self-similar fractals, by understanding the fractal as a boundary of a resistance network.
\end{remark}

%% file: examples.tex

\section{Examples}
\label{sec:examples}

In this section, we introduce the most basic family of examples that illustrate our technical results and exhibit the properties (and support the types of functions) that we have discussed above. 

\begin{exm}[Geometric integer model]\label{def:geometric-integers}
  For a fixed constant $c>1$, let $(\bZ,c^n)$ denote the network with integers for vertices, and with geometrically increasing conductances defined by $\cond_{n-1,n} = c^{\max\{|n|,|n-1|\}}$ so that the network under consideration is
  \linenopax
  \begin{align*}
    \xymatrix{
      \dots \ar@{-}[r]^{c^3}
      & -2 \ar@{-}[r]^{c^2} 
      & -1 \ar@{-}[r]^{c} 
      & 0 \ar@{-}[r]^{c} 
      & 1 \ar@{-}[r]^{c^2} 
      & 2 \ar@{-}[r]^{c^3} 
      & 3 \ar@{-}[r]^{c^4} 
      & \dots
    }
  \end{align*} 
  Fix $o=0$.   
  On this network, the energy kernel is given by
  \linenopax
  \begin{align*}
    v_n(k) = 
    \begin{cases}
      0, &k \leq 0, \\
      \frac{1-r^{k+1}}{1-r}, &1 \leq k \leq n, \\
      \frac{1-r^{n+1}}{1-r}, &k \geq n,
    \end{cases}
    n > 0,
  \end{align*}
  and similarly for $n < 0$.
  Furthermore, the function 
  \linenopax
  \begin{align}\label{eqn:monopole-on-geometric-integers}
    w_o(n) = ar^{|n|}, \q a:= \frac{r}{2(1-r)}
  \end{align}
  defines a monopole, and $h(n) = \operatorname{sgn}(n) (1-w_o(n))$ defines an element of \Harm.
\end{exm}
  
\begin{exm}[Geometric half-integer model]\label{def:geometric-half-integers}
  It is also interesting to consider $(\bZ_+,c^n)$, as this network supports a monopole, but has $\Harm = 0$. 
  \linenopax
  \begin{align*}
    \xymatrix{
      0 \ar@{-}[r]^{c} 
      & 1 \ar@{-}[r]^{c^2} 
      & 2 \ar@{-}[r]^{c^3} 
      & 3 \ar@{-}[r]^{c^4} 
      & \dots
    }
  \end{align*} 
  The monopole can be obtained by rescaling \eqref{eqn:monopole-on-geometric-integers}; just take $a:= \frac{r}{(1-r)}$. There cannot be any nontrivial harmonic functions on this network by \cite[Lem.~5.5]{DGG}, which states that if $h \in \Harm\less\{0\}$, then $h$ has at least two different limiting values at \iy. That is, there exist infinite paths $\cpath_1 = (x_1,x_2,\dots)$ and $\cpath_2 = (y_1,y_2,\dots)$ with $\lim_{j \to \iy} h(x_j) \neq \lim_{j \to \iy} h(y_j)$.
  
  For $k=2,3,\dots$, the network $(\bZ_+,k^n)$ can be thought of as the ``projection'' of the homogeneous tree of degree $k$ $(\sT_k, \tfrac1k \one)$ under a map which sends $x$ to $n \in \bZ$ iff there are $n$ edges between $x$ and $o$.
\end{exm}

\begin{remark}\label{rem:general-integer-networks}
  One can consider more general integer networks, and in this case, $\Harm \neq 0$ for $(\bZ,\cond)$ iff $\sum \cond_{xy}^{-1} < \iy$. In this case, \Harm is spanned by a single bounded function; details appear in \cite{DGG}. Networks of this form have been discussed elsewhere in the literature (for example, \cite[Ex.~3.12, Ex.~4.9]{Kayano88} and \cite[Ex.~3.1, Ex.~3.2]{Kayano84}), but the authors appear to assume that \Lap is self-adjoint. This is not generally the case when \cond is unbounded; in fact, the Laplacian is \emph{not} self-adjoint for Example~\ref{def:geometric-half-integers} or Example~\ref{def:geometric-integers}; see \cite[\S4.2]{SRAMO} or \cite[\S13.4]{OTERN} for further discussion and the explicit computation of defect vectors.
\end{remark}

\begin{exm}[Star networks]\label{exm:star}
  Let $(\sS_m,c^n)$ be a network constructed by conjoining $m$ copies of $(\bZ_+,c^n)$ by identifying the origins of each; let $o$ be the common origin. 
\end{exm}

  Recall from Theorem~\ref{thm:transience} that the boundary term is nontrivial precisely when $\bd \Graph \neq \es$; the presence of a monopole indicates that $\bd \Graph$ contains at least one point. If $\Harm \neq 0$, then there are at least two boundary points; see \cite[Lem.~5.5]{DGG} and Corollary~\ref{thm:Harm-nonzero-iff-multiple-monopoles}.
  
  Example~\ref{exm:star} shows how to construct a network which has a boundary with cardinality $m$. Note that these boundary points can be distinguished by monopoles, by constructing a monopole which is constant everywhere except on one branch. 
  
\begin{exm}[Networks of integer lattices]\label{exm:networks-of-integer-lattices}
  For $d \geq 3$, let \smash{$\{\bZ^d_{(k)}\}_{k=1}^m$} be a collection of $m$ copies of the $d$-dimensional integer lattice $\bZ^d$ with edges between nearest neighbours, and let $o_k$ denote the origin of $\smash{\bZ^d_{(k)}}$. Let $\bZ_m$ be the Cayley graph of the cyclic group of order $m$, and denote its elements by $\{1,2,\dots,m\}$. Now define $(\bZ^d \circledast \bZ_m,\one)$ by identifying $o_k \in \smash{\bZ^d_{(k)}}$ with $k \in \bZ_m$, thus conjoining all the copies of $\bZ^d$. 
  Since $\bZ^d$ is transient, each copy $\bZ^d_{(k)}$ supports a monopole, and hence \Harm has dimension $m-1$ for this network. This is essentially a variation of Example~\ref{exm:star} where $(\bZ_+,c^n)$ is replaced by $\bZ^d$. Note that this is not the same as the Cayley graph of the wreath product $\bZ^d \wr \bZ_m$, which is instead a Diestel-Leader graph; cf.~\cite{Woess05}.
\end{exm}

\begin{exm}[One-sided infinite ladder network]\label{exm:a,b-ladder}
  Consider two copies of the nearest-neighbour graph on the nonnegative integers $\bZ^+$, one with vertices labelled by $\{x_n\}$, and the other with vertices labelled by $\{y_n\}$. Fix two positive numbers $\ga > 1 > \gb > 0$. In addition to the edges $\cond_{x_n, x_{n-1}} = \ga^n$ and $\cond_{y_n, y_{n-1}} = \ga^n$, we also add ``rungs'' to the ladder by defining $\cond_{x_n,y_n} = \gb^n$:
  \linenopax
  \begin{equation}\label{eqn:exm:one-sided-ladder-model}
    \xymatrix{
    *+[l]{x_0} \ar@{-}[r]^{\ga} \ar@{-}[d]_{1}
      &  x_1 \ar@{-}[r]^{\ga^2} \ar@{-}[d]_{\gb} 
      &  x_2 \ar@{-}[r]^{\ga^3} \ar@{-}[d]_{\gb^2} 
      &  x_3 \ar@{-}[r]^{\ga^4} \ar@{-}[d]_{\gb^3}
      &\dots \ar@{-}[r]^{\ga^n}
      &  x_n \ar@{-}[r]^{\ga^{n+1}} \ar@{-}[d]_{\gb^n} & \dots \\
    *+[l]{y_0}\ar@{-}[r]^{\ga} 
      &  y_1 \ar@{-}[r]^{\ga^2} 
      &  y_2 \ar@{-}[r]^{\ga^3} 
      &  y_3 \ar@{-}[r]^{\ga^4} 
      &\dots \ar@{-}[r]^{\ga^n}
      &  y_n \ar@{-}[r]^{\ga^{n+1}} 
      & \dots
    }
  \end{equation}
  This network was suggested to us by Agelos Georgakopoulos. In \cite{bdG}, we show that this example is a one-ended network with nontrivial \Harm, by explicitly constructing a formula for a harmonic function of finite energy on this network. 
\end{exm}

  \version{}{\marginpar{More examples?}}

\begin{exm}[The reproducing kernel on the tree]
  \label{exm:binary-tree:nontrivial-harmonic}
  \label{exm:binary-tree:reproducing-kernel}
  Let $(\sT,\one)$ be the binary tree network as in the top of Figure~\ref{tree-exhaustions} with constant conductance $\cond = \one$. Figure~\ref{fig:tree-repkernels} depicts the embedded image of a vertex $v_x$, as well as its decomposition in terms of \Fin and \Harm. We have chosen $x$ to be adjacent to the origin $o$; the binary label of this vertex would be $x_1$.

  In Figure~\ref{fig:tree-repkernels}, numbers indicate the value of the function at that vertex; artistic liberties have been taken. If vertices $s$ and $t$ are the same distance from $o$, then $|f_x(s)|=|f_x(t)|$ and similarly for $h_x$. Note that $h_x$ provides an example of a nonconstant harmonic function in \HE. 
  It is easy to see that $\lim_{z \to \pm\iy} h_x(z) = \frac12 \pm \frac12$, whence $h_x$ is bounded. 
  

\begin{figure}
  \centering
  \includegraphics{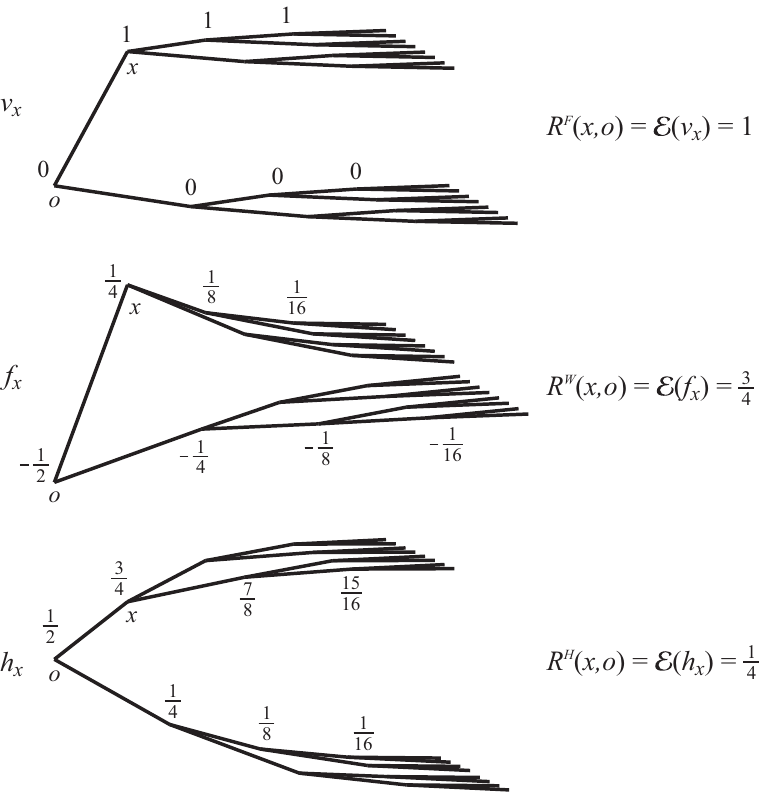}
  \caption{\captionsize The reproducing kernel on the tree with $\cond = \one$. For a vertex $x$ which is adjacent to the origin $o$, this figure illustrates the elements $v_x$, $f_x = \Pfin v_x$, and $h_x = \Phar v_x$; see Example~\ref{exm:binary-tree:nontrivial-harmonic}. }
  \label{fig:tree-repkernels}
\end{figure}

  We can use $h_x$ of Figure~\ref{fig:tree-repkernels} to describe an infinite forest of mutually orthogonal harmonic functions on the binary tree. Let $z \in \sT$ be represented by a finite binary sequence: the root $o$ corresponds to the empty sequence $\es$, and the two vertices connected to it are $0$ and $1$. The neighbours of $0$ are $\es$, $00$ and $01$; the neighbours of $01$ are $0$, $010$, and $011$, etc. Define a mapping $\gf_z: \sT \to \sT $ by prepending, i.e., $\gf_z(x) = zx$. This has the effect of ``rigidly'' translating the the tree so that the image lies on the subtree with root $z$. Then $h_z := h_x \comp \gf_z$ is harmonic and is supported only on the subtree with root $z$. The supports of $h_{z_1}$ and $h_{z_2}$ intersect if and only if $\Im(\gf_{z_i}) \ci \Im(\gf_{z_j})$. For concreteness, suppose it is $\Im(\gf_{z_1}) \ci \Im(\gf_{z_2})$. If they are equal, it is because $z_1=z_2$ and we don't care. Otherwise, compute the dissipation of the induced currents
  \linenopax
  \begin{align*}
    \la \drp h_{z_1}, \drp h_{z_2} \ra_\diss
    = \tfrac12 \negsp \sum_{(x,y) \in \gf_{z_1}(\edges)} \negsp
      \ohm(x,y) \drp h_{z_1}(x,y), \drp h_{z_2}(x,y).
  \end{align*}
  Note that $\drp h_{z_2}(x,y)$ always has the same sign on the subtree with root $z_1 \neq o$, but $\drp h_{z_1}(x,y)$ appears in the dissipation sum positively signed with the same multiplicity as it appears negatively signed. Consequently, all terms cancel and $0 = \la \drp h_{z_1}, \drp h_{z_2} \ra_\diss = \la h_{z_1}, h_{z_2} \ra_\energy$ shows $h_{z_1} \perp h_{z_2}$.

  This family of harmonic functions can be heuristically described by analogy with Haar wavelets.\footnote{Compare to the ``wavelet basis of eigenfunctions'' discussed in \cite{Kig09} (and \cite{Kozyrev02,PearsonBellissard}). These references were brought to our attention by a reader of \cite{OTERN}.} Consider the boundary of the tree as a copy of the unit interval with $h_x$ as the basic Haar mother wavelet; via the ``shadow'' cast by $\lim_{n \to \pm \iy} h_x(x_n) = \pm 1$ (this can be formalized in terms of cocycles). Then $h_z$ is a Haar wavelet localized to the subinterval of the support of its shadow, etc. Of course, this heuristic is a bit misleading, since the boundary is actually isomorphic to $\{0,1\}^\bN$ with its natural cylinder-set topology.
\end{exm}